\documentclass[11pt, letterpaper]{elsarticle}

\usepackage{amsmath}
\usepackage{amssymb}
\usepackage{graphicx}
\usepackage{multirow}
\usepackage{latexsym}
\usepackage{epic}
\usepackage{url}
\usepackage{soul}
\usepackage{afterpage}
\usepackage{subcaption}
\usepackage[unicode, pdftex]{hyperref}
\usepackage{multicol}
\usepackage{bbm}
\usepackage{setspace}
\usepackage{enumitem}
\usepackage{cleveref}
\usepackage[toc]{appendix}
\usepackage[norelsize,ruled, lined,linesnumbered]{algorithm2e}
\usepackage{tikz}
\usepackage{tikz,fullpage}
\usepackage{pgf}
\usepackage{pgfplots}
\usetikzlibrary{
	pgfplots.fillbetween,
}
\usetikzlibrary{arrows,automata}
\usepackage{tkz-berge}
\usepackage{amsthm}
\usepackage[symbol]{footmisc}
\usepackage[mathscr]{euscript}
\usepackage[left=0.9in,top=0.9in,right=0.9in,bottom=0.9in,nohead]{geometry}

\usepackage{empheq}

\DeclareMathOperator{\argmin}{argmin}

\makeatletter
\newtheorem*{rep@theorem}{\rep@title}
\newcommand{\newreptheorem}[2]{%
\newenvironment{rep#1}[1]{%
 \def\rep@title{#2 \ref{##1}}%
 \begin{rep@theorem}}%
 {\end{rep@theorem}}}
\makeatother

\newreptheorem{example}{Example}
\newtheorem{theorem}{Theorem}

\newtheorem{lemma}{Lemma}

\newtheorem{definition}{Definition}

\usepackage{xcolor}
\usepackage[colorinlistoftodos,prependcaption,textsize=tiny]{todonotes}
\usepackage{xargs}

\usepackage{caption}
\usepackage{subcaption}

\usepackage{float}
\usepackage[section]{placeins}
\usepackage{arydshln}
\usepackage{latexsym}
\usepackage{epic}
\usepackage{amsmath,amsthm,amssymb}
\usepackage{graphicx}
\usepackage{url}
\usepackage{pgfplots}\usetikzlibrary{plotmarks}
\usepackage{soul}
\usepackage{afterpage}
\usepackage{cleveref}
\usepackage{bbm}
\usetikzlibrary{decorations.markings}
\usepackage[left=0.9in,top=0.9in,right=0.9in,bottom=0.9in,nohead]{geometry}
\usepackage{setspace}
\usepackage{appendix}

\DeclareGraphicsExtensions{.pdf,.png,.jpg,.bmp}

\allowdisplaybreaks

\begin{document}

\begin{frontmatter}
\title{ A study of distributionally robust mixed-integer programming with Wasserstein metric: on the value of incomplete data }

\author[label1,label2]{Sergey S.~Ketkov\footnote[2]{Corresponding author. Email: sergei.ketkov@business.uzh.ch; phone: +41 078 301 8521. \\\\
The published version of this manuscript is available at: https://doi.org/10.1016/j.ejor.2023.10.018}}

\address[label1]{Laboratory of Algorithms and Technologies for Networks Analysis, HSE University, \\Rodionova st., 136, Nizhny Novgorod, 603093, Russia}
\address[label2]{Department of Business Administration, University of Zurich, Zurich, 8032, Switzerland}
\begin{abstract}
This study addresses a class of linear mixed-integer programming (MILP) problems that involve uncertainty in the objective function parameters. The parameters are assumed to form a random vector, whose probability distribution can only be observed through a finite training data set. Unlike most of the related studies in the literature, we also consider uncertainty in the underlying data set. The data uncertainty is described by a set of linear constraints for each random sample, and the uncertainty in the distribution (for a fixed realization of data) is defined using a type-1 Wasserstein ball centered at the empirical distribution of the data. The overall problem is formulated as a three-level distributionally robust optimization (DRO) problem. First, we prove that the three-level problem admits a single-level MILP reformulation, if the class of loss functions is restricted to biaffine functions. Secondly, it turns out that for several particular forms of data uncertainty, the outlined problem can be solved reasonably fast by leveraging the nominal MILP problem. Finally, we conduct a computational study, where the out-of-sample performance of our model and computational complexity of the proposed MILP reformulation are explored numerically for several application domains.
\end{abstract}

\begin{keyword}
Uncertainty modelling; Distributionally robust optimization; Mixed-integer programming; Wasserstein metric; Incomplete data
\end{keyword}

\end{frontmatter}
\onehalfspace

\section{Introduction} \label{sec: intro}
Distributionally robust optimization (DRO) is a modeling paradigm, in which uncertain problem parameters are described by a \textit{family} (or an \textit{ambiguity set}) of candidate probability distributions that are consistent with the decision-maker's initial information. Specifically, in a standard one-stage DRO problem the decision-maker aims to optimize the expected value of its objective function (or another measure of risk, if the decision-maker is risk-averse) assuming the worst-case distribution of uncertain problem parameters within the ambiguity set; see, e.g., \cite{Delage2010, Esfahani2018, Wiesemann2014}. 

In this paper we employ a DRO approach to a class of data-driven linear mixed-integer programming (MILP) problems with an uncertainty in the objective function parameters. Formally, we consider a stochastic programming problem of the form:
\begin{equation} \label{stochastic programming problem}
\min_{\mathbf{x} \in X} \mathbb{E}_{\mathbb{Q}^*}\{\ell(\mathbf{x}, \mathbf{c})\},
\end{equation}
where 
$X \subseteq \mathbb{R}^{n_1} \times \mathbb{Z}^{n_2}$ is a linear mixed-integer set of feasible decisions \cite{Conforti2014}, $\mathbf{c} \in \mathbb{R}^{n}$, $n = n_1 + n_2$, is a cost vector related to the uncertain problem parameters and $\ell(\mathbf{x}, \mathbf{c})$ is a given loss function. We assume that the nominal (true) distribution $\mathbb{Q}^*$ of the cost vector $\mathbf{c}$ is not known to the decision-maker a priori and can only be observed through a finite training data set. 

 Due to incomplete knowledge of $\mathbb{Q}^*$, the stochastic programming problem (\ref{stochastic programming problem}) cannot be resolved directly. 
However, following the DRO paradigm one may construct a family of probability distributions, which (\textit{i}) is based on the available training data set and (\textit{ii}) contains the nominal distribution~$\mathbb{Q}^*$ with high probability. For example, several recent approaches to modeling data-driven ambiguity sets exploit a distance metric in the space of probability distributions, which is ``centered'' at the empirical distribution of the data; see, e.g., \cite{Bayraksan2015, Bental2013, Esfahani2018, Gao2022}. As a result, for a given family of distributions $\mathcal{Q}$, a standard one-stage DRO problem can be formulated as follows:
\begin{equation} \label{distributionally robust optimization problem}
	\min_{\mathbf{x} \in X} \max_{\mathbb{Q} \in \mathcal{Q}} \mathbb{E}_{\mathbb{Q}}\{\ell(\mathbf{x}, \mathbf{c})\}.
\end{equation}

At the same time, the application of standard DRO models is rather limited when the associated training data set is not completely known (e.g., due to missing data, noise or a particular structure of historical data). In the current study we formulate a new three-level min-max-max optimization model that attempts to handle both uncertainty in the data-generating distribution $\mathbb{Q}^*$ and in the data set obtained from this distribution. While our model is substantially more complicated than the standard DRO model (\ref{distributionally robust optimization problem}), we explore the value of data uncertainty and show that under some additional assumptions our problem can be reformulated as an MILP problem.

\subsection{Related literature} \label{subsec: related literature}

In this section, we discuss several existing data-driven approaches to solving the stochastic programming problem~(\ref{stochastic programming problem}), where the set of feasible decisions, $X$, is either linear mixed-integer or convex. Following the related data-driven stochastic programming formulations in \cite{Delage2010, Esfahani2018}, we assume that only the objective function in (\ref{stochastic programming problem}) but not the constraints is subject to uncertainty. %with either convex or linear mixed-integer sets $X$ of feasible decisions. 
Finally, we distinguish between DRO formulations that use complete or incomplete/partially observable data.

\textbf{DRO with complete data.}
It can be argued that most of the related stochastic programming literature follows the assumption of a \textit{complete training data set}. In other words, it is assumed that the decision-maker in (\ref{stochastic programming problem}) has access to $K \in \mathbb{Z}_{++}$ independently and identically distributed (i.i.d.) observations of the cost vector $\mathbf{c}$ drawn from the nominal distribution $\mathbb{Q}^*$. Under this assumption, the stochastic programming problem (\ref{stochastic programming problem}) can be addressed, for example, using the framework of \textit{sample average approximation} (SAA) or \textit{distributionally robust optimization} (DRO). 

In the former approach the expected value of the loss function in (\ref{stochastic programming problem}) is approximated by the sample mean and the resulting function is optimized over the set $X$ of feasible decisions~\cite{Kleywegt2002}. In general, SAA methods provide deterministic formulations that enjoy strong asymptotic performance guarantees due to the central limit theorem. However, these methods typically provide a rather poor out-of-sample performance when the sample size is not sufficiently large \cite{Bertsimas2018b}. 

On the other hand, it is argued in \cite{Wiesemann2014} that the DRO approach (which is the main focus of this study) enjoys strong justification from decision-theory, where most decision-makers have a low tolerance towards uncertainty in the nominal distribution \cite{Epstein1999}. Most of the solution approaches to min-max DRO problems of the form (\ref{distributionally robust optimization problem}) employ strong duality results for moment problems \cite{Shapiro2001} to obtain equivalent single-level convex \cite{Delage2010, Wiesemann2014} or MILP \cite{Hanasusanto2016, Shang2018, Wang2020} reformulations. Despite the fact that MILP problems are known to be $NP$-hard in general \cite{Garey2002}, this class of problems is of a particular interest by virtue of existing state-of-the-art MILP solvers like CPLEX \cite{Manual1987} and Gurobi \cite{Gurobi2021}. 

As the first example, we refer to the study by Delage and Ye \cite{Delage2010}, who consider convex optimization problems with random parameters in the objective function. The ambiguity set in \cite{Delage2010} is constructed by leveraging a given convex support and confidence sets for the mean and the covariance matrix of the uncertain problem parameters. Then, it is shown that for several particular forms of the loss function, the resulting DRO problems admit finite convex reformulations. 

Despite the fact that moment-based ambiguity sets are rather standard in the DRO literature; see, e.g., \cite{Cheramin2022, Nie2023, Zhang2017c}, it is demonstrated in \cite{Cheng2016} that using the second-order moment constraints with a linear mixed-integer set of feasible decisions results in a \textit{non-linear} single-level MIP reformulation, which can only be approximated by a sequence of semi-definite programming relaxations. %In this regard, several related studies focus on data-driven DRO problems, which admit linear MIP reformulations. 

At the same time, MILP reformulations can be obtained for particular classes of DRO problems with distance-based ambiguity sets. In this regard, we refer to a study by Esfahani and Kuhn~\cite{Esfahani2018}, where the ambiguity set is formed by a ball in the space of (multivariate and non-discrete) probability distributions with respect to the Wasserstein metric. The center of the ball is at the uniform distribution on the training samples, and the radius can be viewed as a decreasing function of the sample size. The authors demonstrate that, if the Wasserstein metric is defined in terms of $l_1$-norm or $l_\infty$-norm, then under some assumptions about the support of $\mathbb{Q}^*$ and the loss function $\ell(\mathbf{x}, \mathbf{c})$, the associated worst-case expectation problem in (\ref{distributionally robust optimization problem}) admits an equivalent dual linear programming reformulation. Therefore, whenever the set of feasible decisions, $X$, is linear mixed-integer, the associated DRO problem can be recast as an MILP problem; see, e.g.,~\cite{Wang2020}.

\textbf{DRO with incomplete/partially observable data.}
First, we refer to the study by Bertsimas~et~al.~\cite{Bertsimas2018}, who analyze a data-driven \textit{robust optimization} approach. That is, instead of the stochastic programming formulation (\ref{stochastic programming problem}), the authors in \cite{Bertsimas2018} consider an associated robust formulation of the form:
\begin{subequations} \label{robust problem}
\begin{align} 
\min_{\mathbf{x} \in X}\max_{\mathbf{c} \in \mathcal{C}} \ell(\mathbf{x}, \mathbf{c}) = & \min_{\mathbf{x}, t} t \\
\mbox{s.t. } & \max_{\mathbf{c} \in \mathcal{C}} \Big( \ell(\mathbf{x}, \mathbf{c}) - t \Big) \leq 0 \label{cons: robust max} \\
& \mathbf{x} \in X,
\end{align}
\end{subequations}
where the loss function $\ell(\mathbf{x}, \mathbf{c})$ is concave in $\mathbf{c}$ and $\mathcal{C}$ is
a convex and compact uncertainty set constructed from $K$ i.i.d. observations of the cost vector $\mathbf{c}$ (according to the nominal distribution $\mathbb{Q}^*$). Then, several hypothesis testing procedures are proposed to construct uncertainty sets with the following two properties: (\textit{i}) the maximization problem in (\ref{cons: robust max}) is tractable; (\textit{ii}) the constraint (\ref{cons: robust max}) is violated with a sufficiently small probability under $\mathbb{Q}^*$. 

Most of the proposed uncertainty sets in \cite{Bertsimas2018} are second-order cone representable or polyhedral with one additional relative entropy constraint and, thus, following our discussion above, their application to discrete optimization problems is rather limited. However, the authors in \cite{Bertsimas2018} consider sampling of data from marginal distributions asynchronously, which may also account \textit{missing data}. In this case the uncertainty set $\mathcal{C}$ is described by box constraints and the overall problem (\ref{robust problem}) admits an MILP reformulation, if, e.g., the loss function $\ell(\mathbf{x}, \mathbf{c})$ is biaffine. 

Next, to the best of our knowledge, there are only a few studies that focus on data-driven DRO formulations of (\ref{stochastic programming problem}) with incomplete data sets. First, we refer to an unpublished preprint by Bennouna and Van Parys \cite{Bennouna2022}, who design distributionally robust formulations that may simultaneously address statistical error, noise and misspecification in the data. In particular, the statistical error is caused by the finite sample size and indicates that we cannot solve the stochastic programming problem (\ref{stochastic programming problem}) exactly. On the other hand, data noise and misspecification are referred to potential measurement errors and corruption of a hopefully small amount of all data, respectively. 
 
The authors in \cite{Bennouna2022} propose a specified DRO approach based on the Kulback-Leibler divergence and the Levy-Prokhorov metric, which is also \textit{robust} in the sense that it may protect against bounded noise and a given fraction of samples that are misspecified. Under some reasonable assumptions about the loss function in (\ref{stochastic programming problem}), the proposed DRO problem is shown to admit a finite convex (but not linear) reformulation for a fixed decision $\mathbf{x} \in X$. Despite the fact that Bennouna and Van Parys \cite{Bennouna2022} do not make any assumptions concerning the set $X$ of feasible decisions, their proposed dual reformulation (Theorem 3.6 in \cite{Bennouna2022}) is only discussed in the context of convex optimization problems, e.g., linear classification and regression. In view of our discussion above, the application of the model in \cite{Bennouna2022} to discrete optimization problems would result in a non-linear MIP problem, which requires more advanced solution techniques than those considered in \cite{Bennouna2022}. %It is probably more important that the model in \cite{Bennouna2022} does not allow to consider some specific forms of noise and misspecification that are dictated by the structure of the underlying optimization problem; we discuss this issue in more detail within the next section.

Finally, we refer to Ren and Bidkhori \cite{Ren2023} for a new data-driven DRO approach, which may handle missing at random (MAR) data. The idea is to use some standard data-driven ambiguity sets based on a distance metric in the space of probability distributions \cite{Bayraksan2015, Esfahani2018}, but replace the empirical distribution of the data by its maximum likelihood estimate (MLE). It is shown in \cite{Ren2023} that, if the data-generating distribution has a \textit{finite discrete support} and $l_1$-norm is utilized to estimate the distance between two distributions, then the problem of finding the MLE estimate reduces to a finite convex optimization problem. Furthermore, by leveraging the obtained estimate of the empirical distribution, standard reformulation techniques can be applied to the resulting DRO problems. Ren and Bidkhori \cite{Ren2023} demonstrate that their approach obeys both asymptotic and finite sample performance guarantees, and consistently outperforms the data imputation approach \cite{Efron1994, Stekhoven2012} for several applied data-driven optimization problems. 

\subsection{Our approach and contributions.} \label{subsec: our approach and contributions}
Perhaps, the major limitation of the proposed data-driven DRO approaches in \cite{Bennouna2022, Bertsimas2018, Ren2023} is that they do not allow to consider any \textit{specific forms} of incomplete data that are dictated by the structure of the underlying optimization problem. More precisely, the studies in \cite{Bertsimas2018, Ren2023} consider only missing data that can be viewed as a particular form of component-wise misspecification. On the other hand, Bennouna and Van Parys \cite{Bennouna2022} introduce bounded noise and misspecification for complete samples drawn from the true joint distribution of the cost vector $\mathbf{c}$. Hence, the modeling approach in \cite{Bennouna2022} is not applicable in cases where noise/misspecification are restricted to a particular subset of components~of~$\mathbf{c}$. 

In order to address the aforementioned gap, in this paper we focus on a class of MILP problems and, especially, on those with a well-defined combinatorial structure. For this class of problems
we propose a new approach to modeling data uncertainty, which allows to consider more flexible data sets tailored to the optimization problem's structure and the data collection process. In particular, our approach is motivated by a number of online combinatorial optimization problem settings, where the decision-maker collects historical data by observing limited information feedback based on its own decisions; see, e.g., \cite{Audibert2014, Bubeck2011}.

\textbf{Approach to data uncertainty.} 
First, we assume that the support set of the cost vector $\mathbf{c}$ is given by a nonempty bounded polytope, i.e.,
\begin{equation} \label{eq: support constraints}
\mathbf{c} \in \mathcal{S}_0 := \Big\{\mathbf{c}' \in \mathbb{R}^{n}: \mathbf{B}^{\mbox{\tiny \upshape (0)}}\mathbf{c}' \leq \mathbf{b}^{\mbox{\tiny \upshape (0)}}\Big\}
\end{equation}
with $\mathbf{B}^{\mbox{\tiny \upshape (0)}} \in \mathbb{R}^{w_0 \times n}$ and $\mathbf{b}^{\mbox{\tiny \upshape (0)}} \in \mathbb{R}^{w_0}$ for some $w_0 \in \mathbb{Z}_{++}$. 
In general, it can be argued that polyhedral support sets are rather standard in DRO; see, e.g., the studies in \cite{Hanasusanto2016, Wiesemann2014}. Secondly, we assume that the nominal distribution $\mathbb{Q}^*$ of the cost vector $\mathbf{c}$ is observed through a finite i.i.d. training data~set 
\begin{equation} \label{eq: data set}
\widehat{\mathbf{C}} = \Big\{ \hat{\mathbf{c}}^{\mbox{\tiny (\itshape k\upshape)}}, \; k \in \mathcal{K}:= \{1, \ldots, K\} \Big\},
\end{equation}
where each sample $\hat{\mathbf{c}}^{\mbox{\tiny (\itshape k\upshape)}} = (\hat{c}^{\mbox{\tiny (\itshape k\upshape)}}_1, \ldots, \hat{c}^{\mbox{\tiny (\itshape k\upshape)}}_n)^\top$, $k \in \mathcal{K}$, is also subject to \textit{linear constraints} of the form:
\begin{equation} \label{eq: data uncertainty}
\hat{\mathbf{c}}^{\mbox{\tiny (\itshape k\upshape)}} \in \mathcal{S}_k := \Big\{\mathbf{c}' \in \mathbb{R}^{n}: \mathbf{B}^{\mbox{\tiny (\itshape k\upshape)}} \mathbf{c}' \leq \mathbf{b}^{\mbox{\tiny (\itshape k\upshape)}}\Big\} \subseteq \mathcal{S}_0
\end{equation}
with $\mathbf{B}^{\mbox{\tiny (\itshape k\upshape)}} \in \mathbb{R}^{w_k \times n}$ and $\mathbf{b}^{\mbox{\tiny (\itshape k\upshape)}} \in \mathbb{R}^{w_k}$, $w_k \in \mathbb{Z}_{++}$. Put differently, by combining the linear data constraints~(\ref{eq: data uncertainty}) with the initial support constraints (\ref{eq: support constraints}), we obtain a subpolytope of $\mathcal{S}_0$, which is simply a singleton when the data set is complete. 

In fact, linear data constraints of the form (\ref{eq: data uncertainty}) can model a rather wide range of incomplete data sets. Among all these data sets we consider the following three major classes of practical interest 
(let $\mathcal{K} := \{1, \ldots, K\}$ and $\mathcal{A} := \{1, \ldots, n\}$ for simplicity of exposition): 
\begin{itemize}
	\item \textbf{Interval uncertainty.} 
	Interval constraints of the form $\hat{c}^{\mbox{\tiny (\itshape k\upshape)}}_{a} \in [l^{\mbox{\tiny (\itshape k\upshape)}}_{a}, u^{\mbox{\tiny (\itshape k\upshape)}}_{a}]$ applied to each $a \in \mathcal{A}$ and $k \in \mathcal{K}$ enable the modeling of bounded noise (we assume that in this case the support constraints (\ref{eq: support constraints}) are also component-wise interval). In contrast to \cite{Bennouna2022}, the magnitude and parameters of the noise depend on the component index $a \in \mathcal{A}$ and can be determined by the decision-maker, either deterministically or randomly.
	\item \textbf{Semi-bandit feedback.}
	To model a situation where some components of $\hat{\mathbf{c}}^{\mbox{\tiny (\itshape k\upshape)}}$ cannot be directly observed by the decision-maker, we combine linear support constraints (\ref{eq: support constraints}) with equality constraints $\hat{c}^{\mbox{\tiny (\itshape k\upshape)}}_{a} = \tilde{c}^{\mbox{\tiny (\itshape k\upshape)}}_{a}$ for some known $\tilde{c}^{\mbox{\tiny (\itshape k\upshape)}}_{a} \in \mathbb{R}_{+}$ and $a \in \mathcal{A}' \subset \mathcal{A}$ (note that any linear equality constraint can always be expressed as two linear inequality constraints). The remaining components are said to be \textit{misspecified}. This type of information feedback is also referred to as ``semi-bandit feedback'' in online learning problem settings \cite{Bubeck2011, Kveton2015}.
	
	\item \textbf{Bandit feedback.} To model a situation where the decision-maker can only observe the total cost with respect to a subset of components $\mathcal{A}' \subset \mathcal{A}$, we introduce linear equality constraints of the form $\sum_{a \in \mathcal{A}'} \hat{c}^{\mbox{\tiny (\itshape k\upshape)}}_{a} = S^{\mbox{\tiny (\itshape k\upshape)}}$ for some $S^{\mbox{\tiny (\itshape k\upshape)}} \in \mathbb{R}_{+}$ and $k \in \mathcal{K}$. %, can be utilized to model a \textit{linear functional dependence} between the components of~$\mathbf{c}$.
 This type of feedback is also known as ``bandit feedback'' in online learning problem settings~\cite{Bubeck2011, Bubeck2012}.
\end{itemize}

 It is clear that by leveraging interval uncertainty and semi-bandit feedback, we can effectively capture both bounded noise and misspecification in our model. Furthermore, the use of semi-bandit and bandit feedback scenarios is strongly justified in the context of \textit{online linear combinatorial optimization problems} (COPs) involving a binary set of feasible decisions $X \subseteq \{0, 1\}^n$. In this setting, the decision-maker may collect historical data by implementing a decision $\mathbf{x} \in X$ and observing either the cost of each nonzero element in $\mathbf{x}$ (semi-bandit feedback) or its total cost (bandit feedback). In addition, we refer to the studies in \cite{Chen2013, Talebi2017, Wen2017} for the application of semi-bandit and bandit feedback scenarios in online COPs such as online advertising, path planning and viral~marketing. These observations provide some intuition behind our choice of $X$ as a linear mixed-integer set. 

Next, in order to handle both uncertainty in the data set and uncertainty in the distribution of the cost vector $\mathbf{c}$, we formulate a three-level min-max-max optimization problem, which is both \textit{robust} (in terms of the data uncertainty) and \textit{distributionally robust} (in terms of the distributional uncertainty). Formally, in the proposed three-level formulation the decision-maker aims to minimize its expected loss by assuming the worst-case possible realization of data with respect to the linear constraints (\ref{eq: data uncertainty}) %from the polyhedral uncertainty sets $\mathcal{S}_k$, $k \in \mathcal{K}$, given by equation (\ref{eq: data uncertainty}) 
and the worst-case realization of the distribution of $\mathbf{c}$ from a predefined ambiguity set of probability distributions. In particular, assuming the worst-case realization of uncertainty aligns with most of the robust and distributionally robust optimization models in the literature; see, e.g., \cite{Bental2002, Esfahani2018, Wiesemann2014}. 

Regarding the ambiguity set, for a fixed realization of data, we focus on a Wasserstein ball w.r.t. $l_1$-norm, centered at the empirical distribution of the data. This choice of norm is consistent with several existing mixed-integer DRO formulations in the literature; see, e.g., \cite{Ji2021, Xie2022}. As a byproduct, some of our results can be slightly modified to capture Wasserstein balls w.r.t. $l_{\infty}$ norm; see Section \ref{sec: solution approach} for further details. Finally, we recall that for both of the outlined ambiguity sets a standard DRO problem of the form (\ref{distributionally robust optimization problem}) with a linear mixed-integer set of feasible decisions $X$ may admit an MILP reformulation~\cite{Esfahani2018, Wang2020}. 

%At the same time, the constraints (\ref{eq: data uncertainty}) are deterministic and, therefore, unlike Ren and Bidkhori \cite{Ren2023}, we do not provide any asymptotic performance guarantees as the sample size tends to infinity. However, in Section \ref{sec: comp study} we illustrate numerically that the out-of-sample performance of our model can be adjusted by an appropriate choice of the polytopes $\mathcal{S}_k$, $k \in \mathcal{K}$. 

\textbf{Contributions.}
The key theoretical result of this study indicates that the proposed three-level problem also admits an MILP reformulation, provided that the class of loss functions in (\ref{stochastic programming problem}) is restricted to \textit{biaffine functions}. Despite the fact that functions of this form are usually not of a particular interest in DRO, they are shown to possess a number of attractive theoretical properties in the context of robust optimization problems \cite{Ben2006, Bertsimas2018}. The linear in $\mathbf{x}$ term of a biaffine loss function may also describe some additional deterministic costs arising, for example, in the context of network interdiction problems~\cite{Israeli2002}. Finally, of even greater importance, \textit{bilinear} loss functions are extensively used as performance measures in online COPs; see, e.g., \cite{Audibert2014}.

Our second contribution is that we consider the particular cases of interval uncertainty, semi-bandit and bandit feedback with a binary set of feasible decisions $X \subseteq \{0, 1\}^n$ and a bilinear loss function $\ell(\mathbf{x}, \mathbf{c}) = \mathbf{c}^\top \mathbf{x}$. It turns out that in this case, under some additional assumptions about the support set (\ref{eq: support constraints}) and the structure of historical data, the three-level problem can be effectively solved by leveraging the nominal MILP problem.

In view of the discussion above, our contributions to the studies in \cite{Bennouna2022, Bertsimas2018, Ren2023} can be summarized as follows:
\begin{itemize}
	\item In contrast to \cite{Bennouna2022, Bertsimas2018, Ren2023}, we propose a new approach to data uncertainty that receives justification from the online learning literature and allows to handle various different forms of historical data that arise in data-driven MILP problems. 
	
	\item In contrast to Bennouona and Van Parys \cite{Bennouna2022}, our three-level optimization model, in general, can be applied to both linear-mixed integer and convex sets of feasible decisions and admits an MILP reformulation in the former case. 
	\item In contrast to Ren and Bidkhori \cite{Ren2023}, who consider distributions with a finite discrete support, we resolve ties to a class of continuous polyhedral support sets that may also capture a linear functional dependence among the components of the cost vector $\mathbf{c}$.
	\item Unlike the studies in \cite{Bennouna2022, Ren2023}, we model the data uncertainty via prespecified sample-wise constraints and, therefore, we do not provide any asymptotic performance guarantees as the sample size tends to infinity. However, it is demonstrated numerically that the out-of-sample performance of our model can be adjusted by an appropriate choice of the linear data constraints~(\ref{eq: data uncertainty}). 
\end{itemize}

 The remainder of the paper is organized as follows.
In Sections \ref{sec: problem} and \ref{subsec: general case}, we formulate the three-level optimization problem and provide its MILP reformulation, respectively. The latter is based on a dual reformulation of the one-stage DRO problem from \cite{Esfahani2018} and Sion's min-max theorem~\cite{Sion1958}. Additionally, we refine the MILP reformulation for the cases of interval uncertainty (semi-bandit feedback) and bandit feedback in Sections \ref{subsec: interval uncertainty} and \ref{subsec: bandit feedback}, respectively. Finally, in Section \ref{sec: comp study} the three-level optimization model is explored numerically for several classes of stochastic combinatorial optimization problems. In particular, we analyze how the form of data uncertainty affects the out-of-sample performance and the computational complexity of our model. 

\textbf{Notation.} All vectors and matrices are labelled by bold letters. A vector of all ones is referred to as~$\mathbf{1}$. We also use subscripts $_+$ and $_{++}$ to define the sets of nonnegative and positive numbers, respectively. A set $\mathcal{K} := \{1, \ldots, K\}$ always refers to the indices of samples in the training data set and a set $\mathcal{A} := \{1, \ldots, n\}$ refers to the indices of components of the cost vector~$\mathbf{c}$. Finally, we denote by $\mathcal{Q}_0(\mathcal{S})$ the space of probability distributions supported on $\mathcal{S} \subseteq \mathbb{R}^{k}$ for some $k \in \mathbb{Z}_{++}$.

\section{Problem formulation} \label{sec: problem}
 Formalizing our discussion in Section \ref{sec: intro}, we consider the stochastic programming problem~(\ref{stochastic programming problem}) under the following additional assumptions: 
 \begin{itemize}
 	\item[\textbf{A1.}] The set of feasible decisions $X$ is linear mixed-integer, i.e., 
 	\begin{equation} 
 	X := \Big\{ (\mathbf{x}^{(1)}, \mathbf{x}^{(2)}) \in \mathbb{R}^{n_1}_+ \times \mathbb{Z}^{n_2}_+: \mathbf{G}_1\mathbf{x}^{(1)} + \mathbf{G}_2 \mathbf{x}^{(2)} \leq \mathbf{g} \Big\}
 	\end{equation} 
 	with $n_1 \in \mathbb{Z}_+$, $n_2 \in \mathbb{Z}_{++}$, $\mathbf{G}_1 \in \mathbb{R}^{m \times n_1}$, $\mathbf{G}_2 \in \mathbb{R}^{m \times n_2}$ and $\mathbf{g} \in \mathbb{R}^{m}$ for some $m \in \mathbb{Z}_{++}$; see, e.g.,~\cite{Conforti2014}. %(in the remainder of the paper we usually do not specify the set $X$ for brevity). 
 	\item[\textbf{A2.}] The loss function $\ell(\mathbf{x}, \mathbf{c})$ in (\ref{stochastic programming problem}) is biaffine, i.e., 
 	\begin{equation}
 	\ell(\mathbf{x}, \mathbf{c}) = \mathbf{c}^\top \mathbf{T} \mathbf{x} + \mathbf{t}_1^\top \mathbf{x} + \mathbf{t}_2^\top \mathbf{c} + t_0,
 	\end{equation}
 	where $\mathbf{T} = \mathbf{T}^\top \in \mathbb{R}^{n \times n}$, $\mathbf{t}_1, \mathbf{t}_2 \in \mathbb{R}^n$ and $t_0 \in \mathbb{R}$.
 	\item[\textbf{A3.}] The cost vector $\mathbf{c}$ belongs to a bounded polyhedral support set $\mathcal{S}_0$ given by equation (\ref{eq: support constraints}). 
 	The nominal distribution $\mathbb{Q}^* \in \mathcal{Q}_0(\mathcal{S}_0)$ of $\mathbf{c}$ can only be observed through a finite i.i.d. training data set (\ref{eq: data set}), where each sample is subject to the linear data constraints (\ref{eq: data uncertainty}). 
 \end{itemize}
 
A detailed discussion of Assumptions \textbf{A1}-\textbf{A3} is provided in Section \ref{subsec: our approach and contributions}. 
In addition, we introduce the following technical definitions. 
First, for a fixed data set $\widehat{\mathbf{C}}$ satisfying Assumption \textbf{A3}, we define an empirical probability distribution obtained from $\widehat{\mathbf{C}}$ as:
\begin{equation} \label{eq: empirical distribution}
\widehat{\mathbb{Q}}_K (\widehat{\mathbf{C}}) := \frac{1}{K}\sum_{k = 1}^K \delta_{\hat{\mathbf{c}}^{\mbox{\tiny (\itshape k\upshape)}}}, 
\end{equation}
where $\delta_{\hat{\mathbf{c}}^{\mbox{\tiny (\itshape k\upshape)}}}$ is the Dirac point mass at the $k$-th training sample $\hat{\mathbf{c}}^{\mbox{\tiny (\itshape k\upshape)}}$ (in the remainder of the paper, the dependence of $\widehat{\mathbb{Q}}_K$ on $\widehat{\mathbf{C}}$ is sometimes omitted for brevity). Next, according to \cite{Kantorovich1958}, we introduce a definition of the Wasserstein distance between two distributions supported on $\mathcal{S}_0$.
\begin{definition} \label{def: Wasserstein distance}
For any $p \in [1, +\infty)$, a type-1 Wasserstein distance between two probability distributions $\mathbb{Q}$ and $\mathbb{Q}'$ on $\mathcal{S}_0$ with respect to $l_p$-norm is defined as:
\begin{equation} \nonumber
\begin{gathered}
W^p(\mathbb{Q}, \mathbb{Q}'):= \inf_{\pi \in \Pi(\mathbb{Q}, \mathbb{Q}')} \int_{\mathcal{S}_0 \times \mathcal{S}_0} \Vert \mathbf{c} - \mathbf{c}' \Vert_p \pi( \mbox{\upshape d}\mathbf{c}, \mbox{\upshape d}\mathbf{c}'),
\end{gathered}	
\end{equation}
where $\Pi(\mathbb{Q}, \mathbb{Q}')$ is a set of all joint distributions of $\mathbf{c} \in \mathcal{S}_0$ and $\mathbf{c}' \in \mathcal{S}_0$ with marginals $\mathbb{Q}$ and $\mathbb{Q}'$, respectively.
\end{definition} 

Since the stochastic programming problem (\ref{stochastic programming problem}) cannot be resolved directly (due to incomplete knowledge of the nominal distribution $\mathbb{Q}^*$ and the associated training data set $\widehat{\mathbf{C}}$), we approximate its solution by introducing a three-level problem of the form:
\begin{equation} \label{three level formulation} \tag{\textbf{F}}
\min_{\mathbf{x} \in X} \max_{\widehat{\mathbf{C}}} \Big\{ \max_{\mathbb{Q} \in \mathcal{Q}(\widehat{\mathbf{C}})} \mathbb{E}_\mathbb{Q} \{\ell(\mathbf{x}, \mathbf{c})\}: \hat{\mathbf{c}}^{\mbox{\tiny (\itshape k\upshape)}} \in \mathcal{S}_k \;\; \forall k \in \mathcal{K} \Big\}, 
\end{equation}
where 
\begin{equation} \label{ambiguity set}
\begin{gathered}
\mathcal{Q}(\widehat{\mathbf{C}}) := \Big\{\mathbb{Q} \in \mathcal{Q}( \mathcal{S}_0 ): \; W^1(\widehat{\mathbb{Q}}_K, \mathbb{Q}) \leq \varepsilon_K \Big\}. \\
\end{gathered}
\end{equation}
The first maximum in (\ref{three level formulation}) indicates that we seek a solution, which is robust to the uncertainty in the data set $\widehat{\mathbf{C}}$ for any fixed decision $\mathbf{x} \in X$. The second maximum in (\ref{three level formulation}) refers to a maximization over the ambiguity set (\ref{ambiguity set}) for a fixed decision $\mathbf{x} \in X$ and a fixed realization of data~$\widehat{\mathbf{C}}$.

The ambiguity set (\ref{ambiguity set}) is defined as a type-1 Wasserstein ball with a radius of $\varepsilon_K > 0$, centered at the empirical distribution of the data; recall (\ref{eq: empirical distribution}). As outlined in Section \ref{sec: intro}, we focus on the Wasserstein distance w.r.t. $l_1$-norm, setting $p = 1$; however, some of our results are applicable to $l_\infty$-norm with $p = \infty$. The Wasserstein radius, $\varepsilon_K$, is defined as a decreasing function in the number of samples, $K$, i.e., the more data is available to the decision-maker, the better it is possible to identify the actual distribution of the uncertain problem parameters. 

We discuss a particular choice of $\varepsilon_K$ in the computational settings; see Section \ref{subsec: test instances}. In the next section we provide an MILP reformulation of the three-level problem (\ref{three level formulation}) and consider three special cases of data uncertainty, namely, interval uncertainty, semi-bandit and bandit feedback. %To the best of our knowledge, the existing estimates of $\varepsilon_K$ are either defined up to a constant factor \cite{Esfahani2018, Fournier2015, La2022} or can be derived from a variation-based concentration inequality \cite{Gao2022b}. We consider a specific form of $\varepsilon_K$ and a linear MIP reformulation of~(\ref{three level formulation}) within the next section. 

\section{Solution approach} \label{sec: solution approach} 
\subsection{General case} \label{subsec: general case}
The reformulation of (\ref{three level formulation}) consists of the following key steps. In the first step, by leveraging duality theory for moment problems \cite{Shapiro2001} we reformulate the third-level maximization problem in (\ref{three level formulation}) as a linear programming problem. The results of this step are based on Assumption \textbf{A2} and the related reformulation of Esfahani and Kuhn \cite{Esfahani2018}. In the second step, we use a version of Sion's min-max theorem \cite{Sion1958} and apply standard linear programming duality to derive an MILP reformulation of the three-level problem (\ref{three level formulation}). The following result holds for the first step.

\begin{lemma} \label{lemma 1}
Assume that a decision $\mathbf{x} \in X$ and a data set $\widehat{\mathbf{C}}$ satisfying Assumption \textbf{A3} are fixed. If, in addition, Assumption \textbf{A2} holds, then the worst-case expectation problem	
\begin{equation} \label{third level problem}
\max_{\mathbb{Q} \in \mathcal{Q}(\widehat{\mathbf{C}})} \mathbb{E}_{\mathbb{Q}} \{\ell(\mathbf{x}, \mathbf{c})\} 
\end{equation} 
admits an equivalent linear programming reformulation of the form:
\begin{subequations} \label{third level LP}
	\begin{align} 
	& \min_{\lambda, \boldsymbol{\nu}} \lambda \varepsilon_K + \frac{1}{K} \sum_{k = 1}^K \Big( \mathbf{b}^{\mbox{\tiny \upshape (0)}\top} \boldsymbol{\nu}^{\mbox{\tiny \upshape (\itshape k\upshape)}} + \hat{\mathbf{c}}^{\mbox{\tiny \upshape (\itshape k\upshape)} \top}(\mathbf{Tx} + \mathbf{t}_2 - \mathbf{B}^{\mbox{\tiny \upshape (0)}\top} \boldsymbol{\nu}^{\mbox{\tiny \upshape (\itshape k\upshape)}})\Big) + \mathbf{t}_1^\top \mathbf{x} + t_0 \\
	\mbox{s.t. } 
	& -\lambda \mathbf{1} \leq \mathbf{T} \mathbf{x} + \mathbf{t}_2 - \mathbf{B}^{\mbox{\tiny \upshape (0)}\top} \boldsymbol{\nu}^{\mbox{\tiny \upshape (\itshape k\upshape)}} \leq \lambda \mathbf{1} \quad \forall{k} \in \mathcal{K} \label{cons: third level LP 1} \\
	& \boldsymbol{\nu}^{\mbox{\tiny \upshape (\itshape k\upshape)}} \geq \mathbf{0} \quad \forall{k} \in \mathcal{K} \label{cons: third level LP 2} \\
	& \lambda \geq 0. \label{cons: third level LP 3}
	\end{align}
\end{subequations}
\begin{proof}
First, for fixed $\mathbf{x}$ and $\widehat{\mathbf{C}}$ the optimization problem (\ref{third level problem}) admits the following dual reformulation (see Theorem 4.2 in the study by Esfahani and Kuhn \cite{Esfahani2018}):
	\begin{subequations} \label{third level dual problem}
	\begin{align} 
	& \min_{ \mathbf{s}, \lambda } \Big( \lambda \varepsilon_K + \frac{1}{K} \sum_{k = 1}^K s_k \Big) \label{obj: third level dual problem} \\
	\mbox{s.t. } & \max_{\mathbf{c}^{\mbox{\tiny \upshape (\itshape k\upshape)}} \in \mathcal{S}_0} \Big (\ell(\mathbf{x}, \mathbf{c}^{\mbox{\tiny \upshape (\itshape k\upshape)}}) - \lambda \Vert \mathbf{c}^{\mbox{\tiny \upshape (\itshape k\upshape)}} - \hat{\mathbf{c}}^{\mbox{\tiny \upshape (\itshape k\upshape)}} \Vert_1 \Big) \leq s_k \quad \forall k \in \mathcal{K} \label{cons: third level dual problem} \\
	& \lambda \geq 0.
	\end{align}
\end{subequations}
 In the following, we fix $k \in \mathcal{K}$ and consider the set of constraints (\ref{cons: third level dual problem}). By utilizing equation (\ref{eq: support constraints}) and the non-negativity constraint $\lambda \geq 0$, we express the maximization problem in the left-hand side of (\ref{cons: third level dual problem}) as a linear programming problem of the form:
	\begin{subequations} \label{fourth level primal problem}
		\begin{align} 
		& \max_{\mathbf{c}^{\mbox{\tiny \upshape (\itshape k\upshape)}}, \mathbf{v}^{\mbox{\tiny \upshape (\itshape k\upshape)}}} \Big(\ell(\mathbf{x}, \mathbf{c}^{\mbox{\tiny \upshape (\itshape k\upshape)}}) - \lambda \sum_{a \in \mathcal{A}} v^{\mbox{\tiny \upshape (\itshape k\upshape)}}_a \Big) \\
		\mbox{s.t. } & \mathbf{B}^{\mbox{\tiny \upshape (0)}} \mathbf{c}^{\mbox{\tiny \upshape (\itshape k\upshape)}} \leq \mathbf{b}^{\mbox{\tiny \upshape (0)}} \label{cons: fourth level primal problem 1} \\
		& -\mathbf{v}^{\mbox{\tiny \upshape (\itshape k\upshape)}} \leq \mathbf{c}^{\mbox{\tiny \upshape (\itshape k\upshape)}} - \hat{\mathbf{c}}^{\mbox{\tiny \upshape (\itshape k\upshape)}} \leq \mathbf{v}^{\mbox{\tiny \upshape (\itshape k\upshape)}}, \label{cons: fourth level primal problem 2}
		\end{align}
	\end{subequations}
where $\ell(\mathbf{x}, \mathbf{c}^{\mbox{\tiny \upshape (\itshape k\upshape)}}) = \mathbf{c}^{\mbox{\tiny \upshape (\itshape k\upshape)}\top} \mathbf{T} \mathbf{x} + \mathbf{t}_1^\top \mathbf{x} + \mathbf{t}_2^\top \mathbf{c}^{\mbox{\tiny \upshape (\itshape k\upshape)}} + t_0$ is a linear function of $\mathbf{c}^{\mbox{\tiny \upshape (\itshape k\upshape)}}$ for a fixed $\mathbf{x} \in X$.

Let $w_0$ be the number of linear constraints in $\mathcal{S}_0$; also let $\boldsymbol{\nu}^{\mbox{\tiny \upshape (\itshape k\upshape)}} \in \mathbb{R}^{w_0}$ and $\boldsymbol{\mu}^{\mbox{\tiny \upshape (\itshape k\upshape)}}_1, \boldsymbol{\mu}^{\mbox{\tiny \upshape (\itshape k\upshape)}}_2 \in \mathbb{R}^{n}$ be dual variables corresponding to the constraints (\ref{cons: fourth level primal problem 1}) and (\ref{cons: fourth level primal problem 2}), respectively, for each $k \in \mathcal{K}$. Then, a dual reformulation of (\ref{fourth level primal problem}) can be expressed as:
	\begin{align*} \label{fourth stage dual problem} 
	& \min_{\boldsymbol{\mu}^{\mbox{\tiny \upshape (\itshape k\upshape)}}_1,\boldsymbol{\mu}^{\mbox{\tiny \upshape (\itshape k\upshape)}}_2, \boldsymbol{\nu}^{\mbox{\tiny \upshape (\itshape k\upshape)}}} \Big(\mathbf{b}^{\mbox{\tiny \upshape (0)}\top} \boldsymbol{\nu}^{\mbox{\tiny \upshape (\itshape k\upshape)}} + \hat{\mathbf{c}}^{\mbox{\tiny \upshape (\itshape k\upshape)} \top} (\boldsymbol{\mu}^{\mbox{\tiny \upshape (\itshape k\upshape)}}_2 - \boldsymbol{\mu}^{\mbox{\tiny \upshape (\itshape k\upshape)}}_1) + \mathbf{t}_1^\top \mathbf{x} + t_0 \Big) \nonumber \\
	\mbox{s.t. } & -\mathbf{Tx} - \mathbf{t}_2 + \mathbf{B}^{\mbox{\tiny \upshape (0)}\top} \boldsymbol{\nu}^{\mbox{\tiny \upshape (\itshape k\upshape)}} - \boldsymbol{\mu}^{\mbox{\tiny \upshape (\itshape k\upshape)}}_1 + \boldsymbol{\mu}^{\mbox{\tiny \upshape (\itshape k\upshape)}}_2 = 0 \nonumber\\
	& \lambda \mathbf{1} - \boldsymbol{\mu}^{\mbox{\tiny \upshape (\itshape k\upshape)}}_1 - \boldsymbol{\mu}^{\mbox{\tiny \upshape (\itshape k\upshape)}}_2 = 0 \nonumber\\
	& \boldsymbol{\mu}^{\mbox{\tiny \upshape (\itshape k\upshape)}}_1, \boldsymbol{\mu}^{\mbox{\tiny \upshape (\itshape k\upshape)}}_2 \geq \mathbf{0} \nonumber\\
	& \boldsymbol{\nu}^{\mbox{\tiny \upshape (\itshape k\upshape)}} \geq \mathbf{0}. \nonumber
	\end{align*}
By eliminating $\boldsymbol{\mu}^{\mbox{\tiny \upshape (\itshape k\upshape)}}_1$ and $\boldsymbol{\mu}^{\mbox{\tiny \upshape (\itshape k\upshape)}}_2$ we obtain the following reformulation of (\ref{fourth level primal problem}) for each $k \in \mathcal{K}$:
 	
\begin{subequations} \label{fourth level dual problem 2}
	\begin{align} 
	& \min_{\boldsymbol{\nu}^{\mbox{\tiny \upshape (\itshape k\upshape)}}} \Big(\mathbf{b}^{\mbox{\tiny \upshape (0)}\top} \boldsymbol{\nu}^{\mbox{\tiny \upshape (\itshape k\upshape)}} + \hat{\mathbf{c}}^{\mbox{\tiny \upshape (\itshape k\upshape)} \top}(\mathbf{Tx} + \mathbf{t}_2 - \mathbf{B}^{\mbox{\tiny \upshape (0)}\top} \boldsymbol{\nu}^{\mbox{\tiny \upshape (\itshape k\upshape)}}) + \mathbf{t}_1^\top \mathbf{x} + t_0 \Big) \\
	\mbox{s.t. } & -\lambda \mathbf{1} \leq \mathbf{Tx} + \mathbf{t}_2 - \mathbf{B}^{\mbox{\tiny \upshape (0)}\top} \boldsymbol{\nu}^{\mbox{\tiny \upshape (\itshape k\upshape)}} \leq \lambda \mathbf{1} \\ 
	& \boldsymbol{\nu}^{\mbox{\tiny \upshape (\itshape k\upshape)}} \geq \mathbf{0}.
	\end{align}
\end{subequations}

Finally, by strong duality the optimal objective function value of (\ref{fourth level dual problem 2}) coincides with the left-hand side of (\ref{cons: third level dual problem}). Taking into account the form of (\ref{cons: third level dual problem}), we can omit the ``min'' operator and shift the constraints in (\ref{fourth level dual problem 2}) to the first level, i.e., the third-level problem (\ref{third level problem}) can be expressed as:
\begin{align*} \label{third level LP preliminary}
	& \min_{\lambda, \boldsymbol{\nu}, \mathbf{s}} \lambda \varepsilon_K + \frac{1}{K} \sum_{k = 1}^K s_k \\
	\mbox{s.t. } & \mathbf{b}^{\mbox{\tiny \upshape (0)}\top} \boldsymbol{\nu}^{\mbox{\tiny \upshape (\itshape k\upshape)}} + \hat{\mathbf{c}}^{\mbox{\tiny \upshape (\itshape k\upshape)} \top}(\mathbf{Tx} + \mathbf{t}_2 - \mathbf{B}^{\mbox{\tiny \upshape (0)}\top} \boldsymbol{\nu}^{\mbox{\tiny \upshape (\itshape k\upshape)}}) + \mathbf{t}_1^\top \mathbf{x} + t_0 \leq s_k \quad \forall{k} \in \mathcal{K} \\
	& -\lambda \mathbf{1} \leq \mathbf{T} \mathbf{x} + \mathbf{t}_2 - \mathbf{B}^{\mbox{\tiny \upshape (0)}\top} \boldsymbol{\nu}^{\mbox{\tiny \upshape (\itshape k\upshape)}} \leq \lambda \mathbf{1} \quad \forall{k} \in \mathcal{K} \\
	& \boldsymbol{\nu}^{\mbox{\tiny \upshape (\itshape k\upshape)}} \geq \mathbf{0} \quad \forall{k} \in \mathcal{K} \\
	& \lambda \geq 0. 
\end{align*}
Then, eliminating the variables $s_k$, $k \in \mathcal{K}$, yields the desired result. 
\end{proof}
\end{lemma} 

By leveraging Lemma \ref{lemma 1}, we can rewrite the three-level formulation (\ref{three level formulation}) as follows:
\begin{equation} \label{three level formulation 2} 
\min_{\mathbf{x} \in X} \max_{\widehat{\mathbf{C}} \in \mathcal{D}_2} \Big\{ \min_{(\lambda, \boldsymbol{\nu}) \in \mathcal{D}_3(\mathbf{x})} f(\mathbf{x}, \widehat{\mathbf{C}}, \lambda, \boldsymbol{\nu}) \Big \}, 
\end{equation}
where the feasible sets, $\mathcal{D}_2$ and $\mathcal{D}_3(\mathbf{x})$, and the objective function, $f(\mathbf{x}, \widehat{\mathbf{C}}, \lambda, \boldsymbol{\nu})$, are given by:
\begin{subequations} \label{eq: feasible sets}
	\begin{align} 
	& \mathcal{D}_2: = \Big\{\widehat{\mathbf{C}} = (\hat{\mathbf{c}}^{\mbox{\tiny \upshape (1\upshape)}}, \ldots, \hat{\mathbf{c}}^{\mbox{\tiny \upshape (\itshape K\upshape)}})^\top \in \mathbb{R}^{K \times n}: \hat{\mathbf{c}}^{\mbox{\tiny \upshape (\itshape k\upshape)}} \in \mathcal{S}_k, \; \forall k \in \mathcal{K} \Big\}, \label{eq: feasible set 1}\\
	& \mathcal{D}_3(\mathbf{x}): = \Big\{(\lambda, \boldsymbol{\nu}) \in \mathbb{R}^{1 \times (w_0 \times K)}: \mbox{ the constraints } (\ref{cons: third level LP 1}), (\ref{cons: third level LP 2}) \mbox{ and } (\ref{cons: third level LP 3}) \mbox{ are satisfied} \Big\}, \label{eq: feasible set 2}\\
	& f(\mathbf{x}, \widehat{\mathbf{C}}, \lambda, \boldsymbol{\nu}) := \lambda \varepsilon_K + \frac{1}{K} \sum_{k = 1}^K \Big( \mathbf{b}^{\mbox{\tiny \upshape (0)}\top} \boldsymbol{\nu}^{\mbox{\tiny \upshape (\itshape k\upshape)}} + \hat{\mathbf{c}}^{\mbox{\tiny \upshape (\itshape k\upshape)} \top}(\mathbf{Tx} + \mathbf{t}_2 - \mathbf{B}^{\mbox{\tiny \upshape (0)}\top} \boldsymbol{\nu}^{\mbox{\tiny \upshape (\itshape k\upshape)}})\Big) + \mathbf{t}_1^\top \mathbf{x} + t_0. \label{eq: obective function modified}
	\end{align}
\end{subequations}
In the following, we show the order of the ``max'' and the second ``min'' operators in (\ref{three level formulation 2}) can be reversed. Furthermore, it turns out that the three-level problem (\ref{three level formulation 2}) admits an MILP reformulation. 

\begin{theorem} \label{theorem 1}
Under Assumptions \textbf{A1}-\textbf{A3} the three-level optimization problem (\ref{three level formulation}) can be reformulated as an MILP problem of the form:	
\begin{subequations} 
	\begin{align} 
	& \min_{\mathbf{x}, \lambda, \boldsymbol{\nu}, \boldsymbol{\gamma}} \lambda \varepsilon_K + \frac{1}{K} \sum_{k = 1}^K \mathbf{b}^{\mbox{\tiny \upshape (0)}\top} \boldsymbol{\nu}^{\mbox{\tiny \upshape (\itshape k\upshape)}} + \sum_{k = 1}^K \mathbf{b}^{\mbox{\tiny \upshape (\itshape k\upshape)}\top} \boldsymbol{\gamma}^{\mbox{\tiny \upshape (\itshape k\upshape)}} + \mathbf{t}^\top_1 \mathbf{x} + t_0 \label{mixed-integer programming reformulation} \tag{\textbf{F}$^{(mip)}$}\\
	\mbox{s.t. } & \frac{1}{K}(\mathbf{Tx} + \mathbf{t}_2 - \mathbf{B}^{\mbox{\tiny \upshape (0)}\top} \boldsymbol{\nu}^{\mbox{\tiny \upshape (\itshape k\upshape)}}) = \mathbf{B}^{\mbox{\tiny \upshape (\itshape k\upshape)}\top} \boldsymbol{\gamma}^{\mbox{\tiny \upshape (\itshape k\upshape)}} \quad \forall k \in \mathcal{K} \nonumber \\
	& -\lambda \mathbf{1} \leq \mathbf{T} \mathbf{x} + \mathbf{t}_2 - \mathbf{B}^{\mbox{\tiny \upshape (0)}\top} \boldsymbol{\nu}^{\mbox{\tiny \upshape (\itshape k\upshape)}} \leq \lambda \mathbf{1} \quad \forall{k} \in \mathcal{K} \nonumber\\
	& \boldsymbol{\nu}^{\mbox{\tiny \upshape (\itshape k\upshape)}}, \boldsymbol{\gamma}^{\mbox{\tiny \upshape (\itshape k\upshape)}} \geq \mathbf{0} \quad \forall{k} \in \mathcal{K} \nonumber\\
	& \lambda \geq 0 \nonumber\\
	& \mathbf{x} \in X.\nonumber
	\end{align}
\end{subequations}

\begin{proof}
In the remainder of the proof we fix $\mathbf{x} \in X$. First, we note that the feasible set $\mathcal{D}_2$ given by equation (\ref{eq: feasible set 1}) is convex and compact as a Cartesian product of $K$ bounded polyhedral uncertainty sets $\mathcal{S}_k$, $k \in \mathcal{K}$. Secondly, the feasible set $\mathcal{D}_3(\mathbf{x})$ given by (\ref{eq: feasible set 2}) is convex since it is described by linear constraints with respect to $\lambda$ and $\boldsymbol{\nu}$; however, we cannot guarantee that $\mathcal{D}_3(\mathbf{x})$ is compact. Finally, it turns out that the objective function $f(\mathbf{x}, \widehat{\mathbf{C}}, \lambda, \boldsymbol{\nu})$ is linear in $\hat{\mathbf{c}}^{\mbox{\tiny \upshape (\itshape k\upshape)}}$ for each $k \in \mathcal{K}$, $\lambda$ and $\boldsymbol{\nu}$. 
 As a result, we can apply Sion's min-max theorem \cite{Sion1958} to the second and the third levels of (\ref{three level formulation 2}), that is, 
\begin{multline}
\max_{\widehat{\mathbf{C}} \in \mathcal{D}_2} \Big\{ \min_{(\lambda, \boldsymbol{\nu}) \in \mathcal{D}_3(\mathbf{x}, \widehat{\mathbf{C}})} f(\mathbf{x}, \widehat{\mathbf{C}}, \lambda, \boldsymbol{\nu}) \Big \} = - \min_{\widehat{\mathbf{C}} \in \mathcal{D}_2} \Big\{ \max_{(\lambda, \boldsymbol{\nu}) \in \mathcal{D}_3(\mathbf{x}, \widehat{\mathbf{C}})} -f(\mathbf{x}, \widehat{\mathbf{C}}, \lambda, \boldsymbol{\nu}) \Big \} = \\ - \max_{(\lambda, \boldsymbol{\nu}) \in \mathcal{D}_3(\mathbf{x}, \widehat{\mathbf{C}})} \Big\{ \min_{\widehat{\mathbf{C}} \in \mathcal{D}_2} -f(\mathbf{x}, \widehat{\mathbf{C}}, \lambda, \boldsymbol{\nu}) \Big \} = \min_{(\lambda, \boldsymbol{\nu}) \in \mathcal{D}_3(\mathbf{x}, \widehat{\mathbf{C}})} \Big\{ \max_{\widehat{\mathbf{C}} \in \mathcal{D}_2} f(\mathbf{x}, \widehat{\mathbf{C}}, \lambda, \boldsymbol{\nu}) \Big \}. \nonumber
\end{multline}

 Next, we observe that the maximization problem 
\begin{equation} \label{third level LP modified}
\max_{\widehat{\mathbf{C}} \in \mathcal{D}_2} f(\mathbf{x}, \widehat{\mathbf{C}}, \lambda, \boldsymbol{\nu})
\end{equation}	
can be viewed as a linear programming problem with respect to $\widehat{\mathbf{C}} \in \mathcal{D}_2$; recall Assumption \textbf{A3}. %Therefore, (\ref{third level LP modified}) admits a linear dual reformulation. 
Let $\boldsymbol{\gamma}^{\mbox{\tiny \upshape (\itshape k\upshape)}}$, $k \in \mathcal{K}$, be dual variables corresponding to the primal constraints $\mathbf{B}^{\mbox{\tiny (\itshape k\upshape)}} \hat{\mathbf{c}}^{\mbox{\tiny (\itshape k\upshape)}} \leq \mathbf{b}^{\mbox{\tiny (\itshape k\upshape)}}$. 
Then, an equivalent dual reformulation of (\ref{third level LP modified}) is given by:
\begin{subequations} \label{second level dual LP}
	\begin{align} 
	& \min_{\lambda, \boldsymbol{\nu}, \boldsymbol{\gamma}} \lambda \varepsilon_K + \frac{1}{K} \sum_{k = 1}^K \mathbf{b}^{\mbox{\tiny \upshape (0)}\top}\boldsymbol{\nu}^{\mbox{\tiny \upshape (\itshape k\upshape)}} + \sum_{k = 1}^K \mathbf{b}^{\mbox{\tiny \upshape (\itshape k\upshape)}\top} \boldsymbol{\gamma}^{\mbox{\tiny \upshape (\itshape k\upshape)}} + \mathbf{t}_1^\top \mathbf{x} + t_0 \\
	\mbox{s.t. } & \frac{1}{K}(\mathbf{Tx} + \mathbf{t}_2 - \mathbf{B}^{\mbox{\tiny \upshape (0)}\top} \boldsymbol{\nu}^{\mbox{\tiny \upshape (\itshape k\upshape)}}) = \mathbf{B}^{\mbox{\tiny \upshape (\itshape k\upshape)}\top} \boldsymbol{\gamma}^{\mbox{\tiny \upshape (\itshape k\upshape)}} \quad \forall k \in \mathcal{K} \\
	& \boldsymbol{\gamma}^{\mbox{\tiny \upshape (\itshape k\upshape)}} \geq \mathbf{0}. \quad \forall k \in \mathcal{K}. 
	\end{align}
\end{subequations}
 In particular, strong duality holds since the set $\mathcal{D}_2$ is bounded and non-empty by construction. 
Eventually, combining the ``min'' operators, as well as the related feasible sets results in a single-level reformulation of the form (\ref{mixed-integer programming reformulation}). 
\end{proof}		
\end{theorem} 

%We also argue that the abovementioned theorem cannot be applied to the case of nonlinear loss functions $\ell(\mathbf{x}, \mathbf{c})$, which provides some intuition behind our formulation of Assumption \textbf{A2}. 

Summarizing the proofs of Lemma \ref{lemma 1} and Theorem \ref{theorem 1}, it can be observed that an equivalent dual reformulation of (\ref{third level problem}) can only be obtained, if the loss function $\ell(\mathbf{x}, \mathbf{c})$ is a \textit{point-wise maximum} of finitely many functions $\ell_r(\mathbf{x}, \mathbf{c})$, $r \in \{1, \ldots, R\}$, that are concave in $\mathbf{c}$, i.e., $\ell(\mathbf{x}, \mathbf{c}) = \max_{r \in \{1,\ldots,R\}} \ell_r(\mathbf{c}, \mathbf{x})$; see,~e.g.,~\cite{Esfahani2018}. We force these functions to be biaffine in both $\mathbf{c}$ and $\mathbf{x}$ to provide an MILP reformulation of the overall problem (\ref{three level formulation}). Moreover, it is rather straightforward to verify that, if $R \geq 2$, then $f(\mathbf{x}, \widehat{\mathbf{C}}, \lambda, \boldsymbol{\nu})$ is convex (but not concave) as a function of~ $\hat{\mathbf{c}}^{\mbox{\tiny \upshape (\itshape k\upshape)}}$, $k \in \mathcal{K}$, that contradicts Sion's min-max theorem. These observations provide some intuition behind Assumption \textbf{A2}, where we set $R = 1$. 

 In addition, we note that the proof of Theorem \ref{theorem 1} does not exploit the structure of $X$. Hence, the three-level problem (\ref{three level formulation}) admits a finite convex reformulation for a convex set of feasible decisions~$X$. Finally, it turns out that slight modifications in the proof of Lemma \ref{lemma 1} imply a result similar to Theorem~\ref{theorem 1} for Wasserstein balls w.r.t. $l_\infty$-norm. 

\subsection{The case of interval uncertainty and semi-bandit feedback } \label{subsec: interval uncertainty} 
In this section, we consider \textit{combinatorial optimization problems} (COPs) with a binary set of feasible decisions $X \subseteq \{0,1\}^{n}$ and a bilinear loss function $\ell(\mathbf{x}, \mathbf{c}) = \mathbf{c}^\top \mathbf{x}$. In this setting, the decision-maker aims to minimize the expected cost of its decision, where both the distribution $\mathbb{Q}^*$ of the cost vector $\mathbf{c}$ and the associated data set $\widehat{\mathbf{C}}$ are subject to uncertainty. We additionally assume that the support constraints (\ref{eq: support constraints}) and the linear data constraints (\ref{eq: data uncertainty}) are given by component-wise \textit{interval constraints}. 
In other words, let
\begin{equation} \label{eq: support constraints interval case}
	\mathcal{S}^{(int)}_0 := \Big\{\mathbf{c}' \in \mathbb{R}^{n}: \mathbf{l} \leq \mathbf{c}' \leq \mathbf{u}\Big\},
\end{equation}
where $\mathbf{0} \leq \mathbf{l} \leq \mathbf{u}$ and 
\begin{equation} \label{eq: data constraints interval case}
	\mathcal{S}^{(int)}_k := \Big\{\mathbf{c}' \in \mathbb{R}^{n}: \mathbf{l}^{\mbox{\tiny \upshape (\itshape k\upshape)}} \leq \mathbf{c}' \leq \mathbf{u}^{\mbox{\tiny \upshape (\itshape k\upshape)}} \Big\},
\end{equation}
for every $k \in \mathcal{K}$, where $\mathbf{l} \leq \mathbf{l}^{\mbox{\tiny \upshape (\itshape k\upshape)}} \leq \mathbf{u}^{\mbox{\tiny \upshape (\itshape k\upshape)}} \leq \mathbf{u}$. 

 This model, including the definitions of $\ell(\mathbf{x}, \mathbf{c})$ and $\mathcal{S}^{(int)}_0$, is similar by construction to the general framework of online combinatorial optimization \cite{Chen2013, Gai2012}. In particular, under the assumption of interval support constraints (\ref{eq: support constraints interval case}), we may also handle the case of semi-bandit feedback by reducing the interval lengths in (\ref{eq: data constraints interval case}) to zero for the observed components of $\mathbf{c}$ and using support constraints~(\ref{eq: support constraints interval case}) for the non-observed components. However, we do to provide a specific reformulation of~(\ref{three level formulation}) for the case of a general polyhedral support set (\ref{eq: support constraints}) and semi-bandit feedback. 

The aforementioned restrictions allow to simplify the proofs of Lemma \ref{lemma 1} and Theorem \ref{theorem 1}. More precisely, we demonstrate that the three-level problem (\ref{three level formulation}) with interval uncertainty/semi-bandit feedback inherits the complexity of the underlying deterministic COP. The following results hold. 

\begin{lemma} \label{lemma 2}
	Assume that $\ell(\mathbf{x}, \mathbf{c}) = \mathbf{c}^\top \mathbf{x}$, $X \subseteq \{0,1\}^{n}$ and the support set $\mathcal{S}^{(int)}_0$ is given by equation (\ref{eq: support constraints interval case}). Then, the optimal objective function value of the third-level problem	
	\begin{equation} \label{third level problem interval case} \nonumber
	\max_{\mathbb{Q} \in \mathcal{Q}(\widehat{\mathbf{C}})} \mathbb{E}_{\mathbb{Q}} \{\ell(\mathbf{x}, \mathbf{c})\} 
	\end{equation} 
	is given by $\min\Big\{\frac{1}{K}\sum_{k = 1}^K \hat{\mathbf{c}}^{\mbox{\tiny \upshape (\itshape k\upshape)}\top} \mathbf{x} + \varepsilon_K; \mathbf{u}^\top \mathbf{x} \Big\}$. 
	\begin{proof}
	Following the proof of Lemma \ref{lemma 1}, we can reformulate the third-level problem (\ref{third level problem}) as follows:
		\begin{subequations} \label{third level dual problem interval case}
			\begin{align} 
			& \min_{ \lambda, \boldsymbol{s} } \Big( \lambda \varepsilon_K + \frac{1}{K} \sum_{k = 1}^K s_k \Big) \label{obj: third level dual problem interval case} \\
			\mbox{s.t. } & \max_{\mathbf{c}^{\mbox{\tiny \upshape (\itshape k\upshape)}} \in [\mathbf{l}, \mathbf{u}]}\Big( \mathbf{c}^{\mbox{\tiny \upshape (\itshape k\upshape)}\top} \mathbf{x} - \lambda \Vert \mathbf{c}^{\mbox{\tiny \upshape (\itshape k\upshape)}} - \hat{\mathbf{c}}^{\mbox{\tiny \upshape (\itshape k\upshape)}} \Vert_1 \Big) \leq s_k \quad \forall k \in \mathcal{K} \label{cons: third level dual problem interval case} \\
			& \lambda \geq 0.
			\end{align}
		\end{subequations}
Then, for each $k \in \mathcal{K}$ the maximum in the left-hand side of (\ref{cons: third level dual problem interval case}) can be expressed as:
\begin{equation}
\begin{gathered} 
\max_{\mathbf{c}^{\mbox{\tiny \upshape (\itshape k\upshape)}} \in [\mathbf{l}, \mathbf{u}]}\Big( \mathbf{c}^{\mbox{\tiny \upshape (\itshape k\upshape)}\top} \mathbf{x} - \lambda \Vert \mathbf{c}^{\mbox{\tiny \upshape (\itshape k\upshape)}} - \hat{\mathbf{c}}^{\mbox{\tiny \upshape (\itshape k\upshape)}} \Vert_1 \Big) = \sum_{a \in \mathcal{A}}\max_{c_a^{\mbox{\tiny \upshape (\itshape k\upshape)}} \in [l_a, u_a]} \Big( c_a^{\mbox{\tiny \upshape (\itshape k\upshape)}} x_a - \lambda |c_a^{\mbox{\tiny \upshape (\itshape k\upshape)}} - \hat{c}_a^{\mbox{\tiny \upshape (\itshape k\upshape)}}| \Big) = \\
\sum_{a \in \mathcal{A}}\max\Big\{\hat{c}_a^{\mbox{\tiny \upshape (\itshape k\upshape)}} x_a; u_a x_a - \lambda (u_a - \hat{c}_a^{\mbox{\tiny \upshape (\itshape k\upshape)}}) \Big\}, \nonumber
\end{gathered}
\end{equation}
where
\begin{equation} \nonumber
\max\Big\{\hat{c}_a^{\mbox{\tiny \upshape (\itshape k\upshape)}} x_a; u_a x_a - \lambda (u_a - \hat{c}_a^{\mbox{\tiny \upshape (\itshape k\upshape)}})\Big\} = \begin{cases} \hat{c}_a^{\mbox{\tiny \upshape (\itshape k\upshape)}} x_a, \mbox{ if } x_a \leq \lambda, \\ u_a x_a - \lambda (u_a - \hat{c}_a^{\mbox{\tiny \upshape (\itshape k\upshape)}}), \mbox{ otherwise. } \end{cases}
\end{equation}	
Here, we exploit the fact that $c_a^{\mbox{\tiny \upshape (\itshape k\upshape)}} x_a - \lambda |c_a^{\mbox{\tiny \upshape (\itshape k\upshape)}} - \hat{c}_a^{\mbox{\tiny \upshape (\itshape k\upshape)}}|$ is a piecewise-linear function of $c_a^{\mbox{\tiny \upshape (\itshape k\upshape)}}$ with $\lambda \geq 0$ and $x_a \geq 0$. Therefore, its maximum can be attained either at $c_a^{\mbox{\tiny \upshape (\itshape k\upshape)}} = \hat{c}_a^{\mbox{\tiny \upshape (\itshape k\upshape)}}$ or at $c_a^{\mbox{\tiny \upshape (\itshape k\upshape)}} = u_a$. 

Next, we consider two particular cases. First, assume that $\lambda \geq 1$ in the dual reformulation (\ref{third level dual problem interval case}). Since $x_a \in \{0, 1\}$ for all $a \in \mathcal{A}$, we observe that:
\begin{equation} \nonumber
\max_{\mathbf{c}^{\mbox{\tiny \upshape (\itshape k\upshape)}} \in [\mathbf{l}, \mathbf{u}]}\Big( \mathbf{c}^{\mbox{\tiny \upshape (\itshape k\upshape)}\top} \mathbf{x} - \lambda \Vert \mathbf{c}^{\mbox{\tiny \upshape (\itshape k\upshape)}} - \hat{\mathbf{c}}^{\mbox{\tiny \upshape (\itshape k\upshape)}} \Vert_1 \Big) = \sum_{a \in \mathcal{A}}\hat{c}_a^{\mbox{\tiny \upshape (\itshape k\upshape)}} x_a
\end{equation}
and, hence, (\ref{third level dual problem interval case}) with an additional requirement $\lambda \geq 1$ has a form:

\begin{align*} \label{third level dual problem interval case 2} 
	& \min_{ \lambda, \mathbf{s} } \Big( \lambda \varepsilon_K + \frac{1}{K} \sum_{k = 1}^K s_k \Big) \nonumber \\
	\mbox{s.t. } & \hat{\mathbf{c}}^{\mbox{\tiny \upshape (\itshape k\upshape)}\top} \mathbf{x} \leq s_k \quad \forall k \in \mathcal{K} \nonumber \\
	& \lambda \geq 1. \nonumber	
\end{align*}	
We conclude that the optimal solution and the optimal objective function value of (\ref{third level dual problem interval case}) are given by $\lambda^* = 1$ and $\frac{1}{K}\sum_{k = 1}^K \hat{\mathbf{c}}^{\mbox{\tiny \upshape (\itshape k\upshape)}\top} \mathbf{x} +~\varepsilon_K$, respectively. %This solution can also be treated as a robust \textit{sample average approximation}, which level of conservatism is controlled by the Wasserstein radius $\varepsilon_K$.

In the second alternative case we assume that $0 \leq \lambda \leq 1$. If $x_a = 0$, then $\max\Big\{\hat{c}_a^{\mbox{\tiny \upshape (\itshape k\upshape)}} x_a; u_a x_a - \lambda (u_a - \hat{c}_a^{\mbox{\tiny \upshape (\itshape k\upshape)}})\Big\} = 0$. Consequently, with $\mathcal{A}_{\mathbf{x}} = \Big\{a \in \mathcal{A}: x_a = 1 \Big\}$ we have:
\begin{equation} 
\begin{gathered}
\max_{\mathbf{c}^{\mbox{\tiny \upshape (\itshape k\upshape)}} \in [\mathbf{l}, \mathbf{u}]}\Big( \mathbf{c}^{\mbox{\tiny \upshape (\itshape k\upshape)}\top} \mathbf{x} - \lambda \Vert \mathbf{c}^{\mbox{\tiny \upshape (\itshape k\upshape)}} - \hat{\mathbf{c}}^{\mbox{\tiny \upshape (\itshape k\upshape)}} \Vert_1 \Big) = \sum_{a \in \mathcal{A}}\max\Big\{\hat{c}_a^{\mbox{\tiny \upshape (\itshape k\upshape)}} x_a; u_a x_a - \lambda (u_a - \hat{c}_a^{\mbox{\tiny \upshape (\itshape k\upshape)}})\Big\} = \\ \sum_{a \in \mathcal{A}_{\mathbf{x}}} \Big(u_a - \lambda (u_a - \hat{c}_a^{\mbox{\tiny \upshape (\itshape k\upshape)}})\Big) = \sum_{a \in \mathcal{A}}\Big(u_a - \lambda (u_a - \hat{c}_a^{\mbox{\tiny \upshape (\itshape k\upshape)}})\Big)x_a, \nonumber
\end{gathered}
\end{equation}
 where the last two equalities follow from the fact that $\lambda \leq 1$ and $x_a \in \{0, 1\}$, $a \in \mathcal{A}$. 
In this case, by eliminating $s_k$, $k \in \mathcal{K}$, the dual problem (\ref{third level dual problem interval case}) can be viewed as: 
\begin{subequations} \label{third level dual problem interval case 3}
	\begin{align} 
	& \min_{\boldsymbol \lambda} \Big( \lambda \varepsilon_K + \frac{1}{K} \sum_{k = 1}^K \sum_{a \in \mathcal{A}}\Big(u_a - \lambda (u_a - \hat{c}_a^{\mbox{\tiny \upshape (\itshape k\upshape)}})\Big)x_a \Big) \label{obj: third level dual problem interval case 3} \\
	\mbox{s.t. } 
	& 0 \leq \lambda \leq 1.
	\end{align}
\end{subequations}
The objective function (\ref{obj: third level dual problem interval case 3}) is linear in $\lambda$ and, hence, $\lambda^* \in \{0, 1\}$. The case of $\lambda^* = 1$ is considered within the first case. Finally, if $\lambda^* = 0$, then the optimal objective function value of (\ref{third level dual problem interval case}) is given by $\mathbf{u}^\top \mathbf{x}$. This observation concludes the proof. 
\end{proof}
\end{lemma}
\begin{theorem} \label{theorem 2}
	Assume that the conditions of Lemma \ref{lemma 2} are satisfied and the data constraints $\mathcal{S}^{(int)}_k$, $k \in \mathcal{K}$, are given by (\ref{eq: data constraints interval case}). Then, an optimal objective function value of (\ref{three level formulation}) coincides with the minimal optimal value of the following two linear COPs:
	\begin{subequations}
		\begin{align}
		& \min_{\mathbf{x} \in X} \frac{1}{K}\sum_{k = 1}^K \mathbf{u}^{\mbox{\tiny \upshape (\itshape k\upshape)}\top} \mathbf{x} + \varepsilon_K, \tag{\textbf{F}$^{(int)}_1$} \label{mixed-integer programming reformulation 1 interval case} \\
		& \min_{\mathbf{x} \in X} \mathbf{u}^\top \mathbf{x}. \tag{\textbf{F}$^{(int)}_2$} \label{mixed-integer programming reformulation 2 interval case}
		\end{align}
	\end{subequations}
\begin{proof}
By leveraging Lemma \ref{lemma 2} we observe that the second-level problem in (\ref{three level formulation}) can be expressed as:
\begin{equation} \label{second-level problem interval case} 
\begin{gathered} \allowdisplaybreaks
\max_{ \widehat{\mathbf{C}} } \Big\{ \max_{\mathbb{Q} \in \mathcal{Q}(\widehat{\mathbf{C}})} \mathbb{E}_\mathbb{Q} \{\mathbf{c}^\top \mathbf{x}\}: \hat{\mathbf{c}}^{\mbox{\tiny (\itshape k\upshape)}} \in \mathcal{S}_k \;\; \forall k \in \mathcal{K} \Big\} = \\ \max_{ \widehat{\mathbf{C}} } \Big\{ \min\big\{\frac{1}{K}\sum_{k = 1}^K \hat{\mathbf{c}}^{\mbox{\tiny \upshape (\itshape k\upshape)}\top} \mathbf{x} + \varepsilon_K; \mathbf{u}^\top \mathbf{x} \big\}: \mathbf{l}^{\mbox{\tiny (\itshape k\upshape)}} \leq \hat{\mathbf{c}}^{\mbox{\tiny (\itshape k\upshape)}} \leq \mathbf{u}^{\mbox{\tiny (\itshape k\upshape)}} \;\; \forall k \in \mathcal{K} \Big\} \geq \\
\min\big\{\frac{1}{K}\sum_{k = 1}^K \mathbf{u}^{\mbox{\tiny \upshape (\itshape k\upshape)}\top} \mathbf{x} + \varepsilon_K; \mathbf{u}^\top \mathbf{x} \big\}.
\end{gathered}
\end{equation}	
 Here, the last inequality holds due to the definition of maximum. Next, we show that the inequality in (\ref{second-level problem interval case}) can be replaced by an equality. 

 Indeed, by setting $\mathbf{s} = \frac{1}{K}\sum_{k = 1}^K \hat{\mathbf{c}}^{\mbox{\tiny \upshape (\itshape k\upshape)}}$ and using (\ref{eq: data constraints interval case}), we observe that 
\begin{equation} \nonumber
\underline{\mathbf{s}} := \frac{1}{K}\sum_{k = 1}^K \mathbf{l}^{\mbox{\tiny \upshape (\itshape k\upshape)}} \leq \mathbf{s} \leq \overline{\mathbf{s}} := \frac{1}{K}\sum_{k = 1}^K \mathbf{u}^{\mbox{\tiny \upshape (\itshape k\upshape)}}.
\end{equation}
Hence, 
\begin{equation} \nonumber
\begin{gathered}
\max_{ \widehat{\mathbf{C}} } \Big\{ \min\big\{\frac{1}{K}\sum_{k = 1}^K \hat{\mathbf{c}}^{\mbox{\tiny \upshape (\itshape k\upshape)}\top} \mathbf{x} + \varepsilon_K; \mathbf{u}^\top \mathbf{x} \big\}: \mathbf{l}^{\mbox{\tiny (\itshape k\upshape)}} \leq \hat{\mathbf{c}}^{\mbox{\tiny (\itshape k\upshape)}} \leq \mathbf{u}^{\mbox{\tiny (\itshape k\upshape)}}, \; \forall k \in \mathcal{K} \Big\} = \\ \max_{ \mathbf{s} } \Big\{ \min\big\{\mathbf{s}^\top \mathbf{x} + \varepsilon_K; \mathbf{u}^\top \mathbf{x} \big\}: \underline{\mathbf{s}} \leq \mathbf{s} \leq \overline{\mathbf{s}} \Big\} = \\ \min\big\{\overline{\mathbf{s}}^\top \mathbf{x} + \varepsilon_K; \mathbf{u}^\top \mathbf{x} \big\} = \min\big\{\frac{1}{K}\sum_{k = 1}^K \mathbf{u}^{\mbox{\tiny \upshape (\itshape k\upshape)}\top} \mathbf{x} + \varepsilon_K; \mathbf{u}^\top \mathbf{x} \big\}.
\end{gathered}
\end{equation}
The second equality above can be checked by contradiction, i.e., for any $\widetilde{\mathbf{s}} \leq \overline{\mathbf{s}}$, $\widetilde{\mathbf{s}} \neq \overline{\mathbf{s}}$, the inequality 
\begin{equation} \nonumber
\min\big\{\widetilde{\mathbf{s}}^\top \mathbf{x} + \varepsilon_K; \mathbf{u}^\top \mathbf{x} \big\} > \min\big\{\overline{\mathbf{s}}^\top \mathbf{x} + \varepsilon_K; \mathbf{u}^\top \mathbf{x} \big\}
\end{equation}
cannot be satisfied under the assumption that $\mathbf{x} \geq \mathbf{0}$. 

We conclude that the three level problem (\ref{three level formulation}) can be viewed as the following optimization problem:
\begin{align}
& \min_{ \mathbf{x} } \Big( \min\{\overline{\mathbf{s}}^\top \mathbf{x} + \varepsilon_K; \mathbf{u}^\top \mathbf{x}\} \Big) \label{min-min problem interval case}\\
\mbox{\upshape s.t. } 
& \mathbf{x} \in X. \nonumber
\end{align}
In (\ref{min-min problem interval case}) we simply find a minimum among $2|X|$ elements. The same result can be obtained by minimizing $\overline{\mathbf{s}}^\top \mathbf{x} + \varepsilon_K$ and $\mathbf{u}^\top \mathbf{x}$ separately with respect to $\mathbf{x} \in X$ and selecting the best optimal value. This observation implies the result.
\end{proof}	
\end{theorem}
We conclude that the considered three-level problem (\ref{three level formulation}) with interval uncertainty/semi-bandit feedback can be solved by leveraging two deterministic linear COPs, (\ref{mixed-integer programming reformulation 1 interval case}) and~(\ref{mixed-integer programming reformulation 2 interval case}). The objective function of (\ref{mixed-integer programming reformulation 1 interval case}) can be seen as a robust sample average approximation of the mean in (\ref{stochastic programming problem}), while the objective function of (\ref{mixed-integer programming reformulation 2 interval case}) represents the worst-case expected cost that can be incurred by the decision-maker. Furthermore, if the associated nominal COP is polynomially solvable, then (\ref{three level formulation}) can also be solved in polynomial~time. In conclusion, we note that, in contrast to Theorem \ref{theorem 1}, the proofs of Lemma \ref{lemma 2}, Theorem \ref{theorem 2} and subsequent results in Section \ref{subsec: bandit feedback} cannot be modified in a straightforward way to capture Wasserstein balls w.r.t. $l_\infty$-norm. % with the cost vectors given by $\frac{1}{K}\sum_{k = 1}^K \mathbf{u}^{\mbox{\tiny \upshape (\itshape k\upshape)}}$ and $\mathbf{u}$, respectively. 
\subsection{The case of bandit feedback} \label{subsec: bandit feedback}
In this section we also assume that $X \subseteq \{0, 1\}^{n}$ and $\ell(\mathbf{x}, \mathbf{c}) = \mathbf{c}^\top \mathbf{x}$. However, in order to provide a problem-specific reformulation of (\ref{three level formulation}) for the case of bandit feedback, we need to make the following additional assumptions. First, we suppose that any decision $\mathbf{x} \in X$ has exactly $h \in \mathbb{Z}_{++}$ non-zero components, i.e., $\Vert \mathbf{x} \Vert_1 = h$, and define the support constraints (\ref{eq: support constraints}) as: 
\begin{equation} \label{eq: support constraints bf}
	\mathcal{S}^{(bf)}_0 := \Big\{\mathbf{c}' \in \mathbb{R}^{n}: \mathbf{0} \leq \mathbf{c}' \leq \mathbf{1}\Big\},
\end{equation}
 That is, in contrast to the case of interval uncertainty, each component of the cost vector $\mathbf{c}$ is assumed to belong to the same interval, $[0, 1]$. 

Next, for a given set of decisions $\mathbf{x}^{\mbox{\tiny \upshape (\itshape k\upshape)}}$, $k \in \mathcal{K}$, the linear data constraints (\ref{eq: data uncertainty}) are defined as:
\begin{align}
	\mathcal{S}^{(bf)}_k := \Big\{\mathbf{c}' \in \mathbb{R}^{n}: \mathbf{c}'^\top \mathbf{x}^{\mbox{\tiny \upshape (\itshape k\upshape)}} = S^{\mbox{\tiny \upshape (\itshape k\upshape)}}\Big\} \cap \mathcal{S}^{(bf)}_0, \label{eq: data uncertainty bf 1}
\end{align}
 where $S^{\mbox{\tiny \upshape (\itshape k\upshape)}} \in \mathbb{R}_+$ denotes the total cost associated with $\mathbf{x}^{\mbox{\tiny \upshape (\itshape k\upshape)}}$. Finally, we make the following additional assumption about the set of decisions $\mathbf{x}^{\mbox{\tiny \upshape (\itshape k\upshape)}}$, $k \in \mathcal{K}$ (a pair of decisions $\mathbf{x}' \in X$ and $\mathbf{x}'' \in X$ is called non-overlapping, if $x'_a + x''_a \in \{0, 1\}$ for each~$a \in \mathcal{A}$): 
\begin{itemize}
	\item[\textbf{A4}.] The set of decisions $\mathbf{x}^{\mbox{\tiny \upshape (\itshape k\upshape)}}$, $k \in \mathcal{K}$, is comprised of $V$ mutually non-overlapping decisions $\widetilde{\mathbf{x}}^{\mbox{\tiny \upshape (\itshape v\upshape)}}$, $v \in \{1, \ldots, V\}$, where the total cost of each decision $\widetilde{\mathbf{x}}^{\mbox{\tiny \upshape (\itshape v\upshape)}}$ is observed $n_v \geq 1$ times, $\sum_{v = 1}^V n_v = K$.
\end{itemize}

The aforementioned assumptions (except for Assumption \textbf{A4}) align with those made in the related online combinatorial optimization literature; see, e.g., \cite{Audibert2014}. With respect to Assumption \textbf{A4}, it may have limited applicability to real historical data sets, where the observed decisions, e.g., routes in the network, are supposed to overlap with each other. However, as we will demonstrate later, the reformulation of (\ref{three level formulation}) under Assumption \textbf{A4} can offer valuable practical insights even for the case when Assumption \textbf{A4} is relaxed. The following result holds.

\begin{theorem} \label{theorem 3}
	Let $X \subseteq \{0,1\}^{n}$, $\ell(\mathbf{x}, \mathbf{c}) = \mathbf{c}^\top \mathbf{x}$ and $\Vert \mathbf{x} \Vert_1 = h$ for any $\mathbf{x} \in X$. Assume that the set of decisions $\mathbf{x}^{\mbox{\tiny \upshape (\itshape k\upshape)}}$, $k \in \mathcal{K}$, satisfies Assumption \textbf{A4} and the empirical mean of the cost of $\widetilde{\mathbf{x}}^{\mbox{\tiny \upshape (\itshape v\upshape)}}$, $v \in \{1, \ldots, V\}$, is given by $\overline{S}_v$. Then, the optimal objective function value of (\ref{three level formulation}) subject to (\ref{eq: support constraints bf}) and (\ref{eq: data uncertainty bf 1}) coincides~with 
\begin{equation} \label{three-level reformulation bf} \tag{\textbf{F}$^{(sbf)}$}
\min \Big\{ \frac{1}{K} \min_{v \in \{1, \ldots V\}}\big\{ n_v \overline{S}_v + (K - n_v)h\big \} + \varepsilon_K; h\Big\}
\end{equation}
and an optimal decision $\mathbf{x}^* = \widetilde{\mathbf{x}}^{\mbox{\tiny \upshape (\itshape v*\upshape)}}$, where $v^* \in \argmin_{v \in \{1, \ldots V\}}\big\{ n_v \overline{S}_v + (K - n_v)h\big\}$.
\begin{proof}
	First, we note that the conditions of Lemma \ref{lemma 2} are satisfied and, therefore, the optimal objective function value of the third-level problem	
	\begin{equation} \label{third level problem interval case 2} \nonumber
	\max_{\mathbb{Q} \in \mathcal{Q}(\widehat{\mathbf{C}})} \mathbb{E}_{\mathbb{Q}} \{\ell(\mathbf{x}, \mathbf{c})\} 
	\end{equation} 
	is given by $\min\{\frac{1}{K}\sum_{k = 1}^K \hat{\mathbf{c}}^{\mbox{\tiny \upshape (\itshape k\upshape)}\top} \mathbf{x} + \varepsilon_K; \mathbf{u}^\top \mathbf{x} \} = \min\{\frac{1}{K}\sum_{k = 1}^K \hat{\mathbf{c}}^{\mbox{\tiny \upshape (\itshape k\upshape)}\top} \mathbf{x} + \varepsilon_K; h \}$. 
	Hence, the second-level problem in (\ref{three level formulation}) can be expressed as:
	\begin{subequations}
	\begin{align} \label{eq: second-level problem bf}
	& \max_{ \widehat{\mathbf{C}} } \min_{\mathbf{y} \in Y} \Big\{ \Big(\frac{1}{K}\sum_{k = 1}^K \hat{\mathbf{c}}^{\mbox{\tiny \upshape (\itshape k\upshape)}\top} \mathbf{x} + \varepsilon_K\Big)y_1 + hy_2\Big\} \\ \mbox{ s.t. } & \mathbf{0} \leq \hat{\mathbf{c}}^{\mbox{\tiny \upshape (\itshape k\upshape)}} \leq \mathbf{1} \\ & \hat{\mathbf{c}}^{\mbox{\tiny \upshape (\itshape k\upshape)}\top} \mathbf{x}^{\mbox{\tiny \upshape (\itshape k\upshape)}} = S^{\mbox{\tiny \upshape (\itshape k\upshape)}} \quad \forall k \in \mathcal{K},
	\end{align}
 \end{subequations}
	where $\mathbf{x}^{\mbox{\tiny \upshape (\itshape k\upshape)}} = \widetilde{\mathbf{x}}^{\mbox{\tiny \upshape (\itshape v\upshape)}}$ for some $v \in \{1, \ldots, V\}$ and $Y = \{(y_1, y_2) \in \mathbb{R}^2_{+}: y_1 + y_2 = 1\}$. 
	
	By Sion's min-max theorem \cite{Sion1958}, we can reverse the order of operators in (\ref{eq: second-level problem bf}), which results in the following reformulation:
		\begin{align} \label{eq: second-level problem bf 2} \nonumber
		& \min_{\mathbf{y} \in Y} \Big\{ \Big(\frac{1}{K}\sum_{k = 1}^K \max_{\hat{\mathbf{c}}^{\mbox{\tiny \upshape (\itshape k\upshape)}} \in 	\mathcal{S}^{(bf)}_k} \hat{\mathbf{c}}^{\mbox{\tiny \upshape (\itshape k\upshape)}\top} \mathbf{x} + \varepsilon_K\Big)y_1 + hy_2\Big\}.
		\end{align}
	Next, it is rather easy to show that 
	\begin{equation} \nonumber
	\max_{\hat{\mathbf{c}}^{\mbox{\tiny \upshape (\itshape k\upshape)}} \in 	\mathcal{S}^{(bf)}_k} \hat{\mathbf{c}}^{\mbox{\tiny \upshape (\itshape k\upshape)}\top} \mathbf{x} = \min\Big\{S^{\mbox{\tiny \upshape (\itshape k\upshape)}} + \sum_{a \in \mathcal{A}} \max\{0, x_a - x^{\mbox{\tiny \upshape (\itshape k\upshape)}}_a\}; h\Big\} 
	\end{equation}
and, therefore, the three-level problem (\ref{three level formulation}) admits the following equivalent reformulation:
\begin{equation} \label{eq: first-level problem bf}
\begin{gathered}
\min_{\mathbf{x} \in X}\min_{\mathbf{y} \in Y} \Big\{ \Big(\frac{1}{K}\sum_{k = 1}^K \min\big\{S^{\mbox{\tiny \upshape (\itshape k\upshape)}} + \sum_{a \in \mathcal{A}} \max\{0, x_a - x^{\mbox{\tiny \upshape (\itshape k\upshape)}}_a\}; h \big \} + \varepsilon_K\Big)y_1 + hy_2\Big\} = \\
\min_{\mathbf{y} \in Y} \Big\{ \Big(\frac{1}{K} \min_{\mathbf{x} \in X} \sum_{k = 1}^K \min\big\{S^{\mbox{\tiny \upshape (\itshape k\upshape)}} + \sum_{a \in \mathcal{A}} \max\{0, x_a - x^{\mbox{\tiny \upshape (\itshape k\upshape)}}_a\}; h \big \} + \varepsilon_K\Big)y_1 + hy_2\Big\}.
\end{gathered}
\end{equation}
Here, the last equality follows from the fact that the minimization over $\mathbf{x} \in X$ and $\mathbf{y} \in Y$ can be realized in any predefined order.
 	
Finally, we consider two particular cases. First, let $\mathbf{x} = \widetilde{\mathbf{x}}^{\mbox{\tiny \upshape (\itshape v\upshape)}}$ for some $v \in \{1, \ldots, V\}$. Taking into account that $\widetilde{\mathbf{x}}^{\mbox{\tiny \upshape (\itshape v\upshape)}} \in \{0, 1\}^n$ and $\Vert \widetilde{\mathbf{x}}^{\mbox{\tiny \upshape (\itshape v\upshape)}} \Vert_1 = h$, we conclude that for any $k \in \mathcal{K}$
\begin{equation} \nonumber
\sum_{a \in \mathcal{A}}\max\{0, x_a - x^{\mbox{\tiny \upshape (\itshape k\upshape)}}_a\} = \begin{cases} 0, \mbox{ if } \mathbf{x} = \mathbf{x}^{\mbox{\tiny \upshape (\itshape k\upshape)}}, \\
h, \mbox{ if } \mathbf{x} = \mathbf{x}^{\mbox{\tiny \upshape (\itshape j\upshape)}}, \; j \neq k
\end{cases}
\end{equation}
and, thus, 
\begin{equation} \nonumber
R_v := \sum_{k = 1}^K \min\big\{S^{\mbox{\tiny \upshape (\itshape k\upshape)}} + \sum_{a \in \mathcal{A}} \max\{0, x_a - x^{\mbox{\tiny \upshape (\itshape k\upshape)}}_a\}; h \big \} = n_v \overline{S}_v + (K - n_v)h;
\end{equation}
we recall that $S^{\mbox{\tiny \upshape (\itshape k\upshape)}} \in [0, h]$ by the support constraints (\ref{eq: support constraints bf}).
%Furthermore, $S^{\mbox{\tiny \upshape (\itshape k\upshape)}} \in [0, h]$ by the support constraints (\ref{eq: support constraints bf}) with $\mathbf{l} = 0$ and $\mathbf{u} = \mathbf{1}$. 

Secondly, for $\mathbf{x} \in X$ such that $\mathbf{x} \neq \widetilde{\mathbf{x}}^{\mbox{\tiny \upshape (\itshape v\upshape)}}$, $v \in \{1, \ldots, V\}$, we demonstrate that $\mathbf{x}$ cannot be an optimal solution of (\ref{three level formulation}). In this regard, the total costs $S^{\mbox{\tiny \upshape (\itshape k\upshape)}}$, $k \in \mathcal{K}$, are divided into $V$ distinct groups associated with the decisions $\widetilde{\mathbf{x}}^{\mbox{\tiny \upshape (\itshape v\upshape)}}$, $v \in \{1, \ldots, V\}$. Formally, let $\widetilde{S}^{\mbox{\tiny \upshape (\itshape v, j\upshape)}}$, $j \in \{1, \ldots, n_v\}$, be the $j$-th observation of the total cost of $\widetilde{\mathbf{x}}^{\mbox{\tiny \upshape (\itshape v\upshape)}}$. Then, by setting $\mathcal{A}_{\mathbf{x}} := \{a \in \mathcal{A}: x_a = 1 \}$ we observe that: 
\begin{equation} \label{eq: theorem 3 proof 1} 
\begin{gathered} 
\sum_{k = 1}^K \min\big\{S^{\mbox{\tiny \upshape (\itshape k\upshape)}} + \sum_{a \in \mathcal{A}} \max\{0, x_a - x^{\mbox{\tiny \upshape (\itshape k\upshape)}}_a\}; h \big \} = \sum_{k = 1}^K \min\big\{S^{\mbox{\tiny \upshape (\itshape k\upshape)}} + |\mathcal{A}_{\mathbf{x}} \setminus \mathcal{A}_{\mathbf{x}^{\mbox{\tiny \upshape (\itshape k\upshape)}}}|; h\big \} = \\
\\ \sum_{v = 1}^V \sum_{j = 1}^{n_v} \min\big\{ \widetilde{S}^{\mbox{\tiny \upshape (\itshape v, j\upshape)}} + |\mathcal{A}_{\mathbf{x}} \setminus \mathcal{A}_{\;\widetilde{\mathbf{x}}^{\mbox{\tiny \upshape (\itshape v\upshape)}}}|; h\big \}.
\end{gathered}
\end{equation}

In the following, for each $v \in \{1, \ldots, V\}$ and $j \in \{1, \ldots, n_v\}$ we introduce coefficients $\alpha_{v,j } \in [0, 1]$ and $\beta_v \in [0,1]$ such that $\widetilde{S}^{\mbox{\tiny \upshape (\itshape v, j\upshape)}} = \alpha_{v,j}h$, $|\mathcal{A}_{\mathbf{x}} \cap \mathcal{A}_{\widetilde{\mathbf{x}}^{\mbox{\tiny \upshape (\itshape v\upshape)}}}| = \beta_v h$ and $\sum_{v = 1}^V \beta_v = \beta \leq 1$.
Furthermore, without loss of generality for any fixed $v \in \{1, \ldots, V\}$ there exist an index $r_v \in \{0, \ldots, n_v\}$ such that
\begin{equation} \nonumber
\alpha_{v,0} := 0 \leq \alpha_{v,1} \leq \ldots \leq \alpha_{v,r_v} \leq \beta_v \leq \alpha_{v,r_v + 1} \leq \ldots \leq \alpha_{v,n_v} \leq 1.
\end{equation}
As a result, 
\begin{equation} \label{eq: theorem 3 proof 2}
\begin{gathered} 
\sum_{v = 1}^V \sum_{j = 1}^{n_v} \min\big\{ \widetilde{S}^{\mbox{\tiny \upshape (\itshape v, j\upshape)}} + |\mathcal{A}_{\mathbf{x}} \setminus \mathcal{A}_{\;\widetilde{\mathbf{x}}^{\mbox{\tiny \upshape (\itshape v\upshape)}}}|; h\big \} = \sum_{v = 1}^V \sum_{j = 1}^{n_v}\min\big \{(\alpha_{v,j} + 1 - \beta_v)h ; h\big \} = \\
\sum_{v = 1}^V \Big(\sum_{j = 1}^{r_v} (\alpha_{v,j} + 1 - \beta_v)h + \sum_{j = r_v + 1}^{n_v}h \Big) \geq \sum_{v = 1}^V \Big(\sum_{j = 1}^{r_v} (\alpha_{v,j} + 1 - \widetilde{\beta}_v)h + \sum_{j = r_v + 1}^{n_v}h \Big) = \\
\sum_{v = 1}^V \Big(\sum_{j = 1}^{r_v} \alpha_{v,j}h + \sum_{j = r_v + 1}^{n_v}\widetilde{\beta}_v h + \sum_{j = 1}^{n_v}(1 - \widetilde{\beta}_v)h \Big),
\end{gathered}
\end{equation}
where $\widetilde{\beta}_v = \frac{\beta_v}{\beta} \in [0, 1]$ and $\sum_{v = 1}^V \widetilde{\beta}_v = 1$ (if $\beta = 0$, then $|\mathcal{A}_{\mathbf{x}} \cap \mathcal{A}_{\;\widetilde{\mathbf{x}}^{\mbox{\tiny \upshape (\itshape v\upshape)}}}| = 0$ for any $v \in \{1, \ldots, V\}$ and $\mathbf{x}$ is clearly suboptimal as its worst-case expected cost equals $h$).
In particular, by using the fact that $\alpha_{v,k} \in [0, 1]$ and $\beta_v \in [0,1]$ the two terms in (\ref{eq: theorem 3 proof 2}) can be estimated as follows:
\begin{align} \nonumber
& \sum_{v = 1}^V \Big(\sum_{j = 1}^{r_v} \alpha_{v,j}h + \sum_{j = r_v + 1}^{n_v}\widetilde{\beta}_v h\Big) \geq \sum_{v = 1}^V \Big(\widetilde{\beta}_v \sum_{j = 1}^{r_v} \alpha_{v,j}h + \widetilde{\beta}_v \sum_{j = r_v + 1}^{n_v} \alpha_{v,j} h\Big) =
\sum_{v = 1}^V \widetilde{\beta}_v n_v \overline{S}_v, \\
& \sum_{v = 1}^V \sum_{j = 1}^{n_v}(1 - \widetilde{\beta}_v)h = \Big(K - \sum_{v = 1}^V \widetilde{\beta}_v n_v \Big)h. \nonumber
\end{align}

By combining the outlined observations with (\ref{eq: theorem 3 proof 1}) we conclude that:
\begin{equation} \nonumber
\sum_{k = 1}^K \min\big\{S^{\mbox{\tiny \upshape (\itshape k\upshape)}} + \sum_{a \in A} \max\{0, x_a - x^{\mbox{\tiny \upshape (\itshape k\upshape)}}_a\}; h \big \} \geq \sum_{v = 1}^V \widetilde{\beta}_v \Big(n_v \overline{S}_v + (K - n_v)h\Big) = \sum_{v = 1}^V \widetilde{\beta}_v R_v \geq \min_{v \in \{1, \ldots v\}} R_v,
\end{equation}
where the last inequality with (\ref{eq: first-level problem bf}) yields that $\mathbf{x}$ is a suboptimal solution of (\ref{three level formulation}). Also, from (\ref{eq: first-level problem bf}) we observe that an optimal objective function value of (\ref{three level formulation}) coincides with 
\begin{equation} \nonumber
\min_{\mathbf{y} \in Y}\Big(\frac{1}{K} \min_{v \in \{1, \ldots V\}} R_v + \varepsilon_K \Big) y_1 + hy_2 = \min\Big\{\frac{1}{K} \min_{v \in \{1, \ldots V\}} R_v + \varepsilon_K; h\Big\}.
\end{equation}
This observation implies the result.
\end{proof}
\end{theorem}

Theorem \ref{theorem 3} states the three-level problem (\ref{three level formulation}) with bandit feedback and non-overlapping decisions can always be efficiently solved by enumerating the values of $R_v := n_v \overline{S}_v + (K - n_v)h$, $v \in \{1, \ldots, V\}$.
 Importantly, this result is established with the only assumption that $\Vert \mathbf{x} \Vert_1 =~h$ for any $\mathbf{x} \in X \subseteq \{0, 1\}^n$ and without any additional assumptions on the structure of~$X$. Furthermore, it turns out that the three-level problem~(\ref{three level formulation}) with bandit feedback can be solved rather effectively, even when Assumption~\textbf{A4} is relaxed. That is, in the next section, we demonstrate numerically that the MILP reformulation (\ref{mixed-integer programming reformulation}) for a class of network routing problems with bandit feedback can be well approximated by its LP relaxation. 

\section{Computational study} \label{sec: comp study}
In this section we aim to analyze the practical applicability of the three-level problem (\ref{three level formulation}) by considering the value of data uncertainty and the computational complexity of the MILP reformulation~(\ref{mixed-integer programming reformulation}). To accomplish this, we explore different forms of data uncertainty, including interval uncertainty and semi-bandit/bandit feedback, and apply them to three specific classes of stochastic combinatorial optimization problems: the \textit{sorting problem}, the \textit{shortest path problem}, and the \textit{maximum coverage problem}; see, e.g., \cite{Korte2011}. 

The remainder of this section is organized as follows. In Section \ref{subsec: test instances} we describe our test instances including the classes of combinatorial optimization problems and the form of data-generating distribution. In Section \ref{subsec: results and discussion} the related computational results and their discussion are provided. 
 
\subsection{Test instances and computational settings} \label{subsec: test instances}
\textbf{Measure of performance.} The \textit{out-of-sample performance} of our model for a given nominal distribution $\mathbb{Q}^* \in \mathbb{Q}_0( \mathcal{S}_0 )$ is evaluated by leveraging a \textit{nominal relative loss}. In other words, if $\widetilde{\mathbf{x}}^* \in X$ denotes an optimal solution of (\ref{three level formulation}), we assess its quality using the following equation:
\begin{eqnarray} \label{eq: nominal relative loss}
\rho(\widetilde{\mathbf{x}}^*, \mathbb{Q}^*) = \frac{\mathbb{E}_{\mathbb{Q}^*}\{\ell(\widetilde{\mathbf{x}}^*, \mathbf{c})\}}{\min_{\mathbf{x} \in X} \mathbb{E}_{\mathbb{Q}^*}\{\ell(\mathbf{x}, \mathbf{c})\}}.
\end{eqnarray}
It is worth noting that by design, the nominal relative loss $\rho(\widetilde{\mathbf{x}}^*, \mathbb{Q}^*)$ is greater than or equal to 1.

\textbf{Classes of problems and data uncertainty.}
%We consider two classical polynomially solvable COPs, namely, the \textit{sorting problem} and the \textit{shortest path problem}. 
 Taking into account Theorems \ref{theorem 2} and \ref{theorem 3}, the focus of our numerical study is on COPs with a binary set of feasible decisions $X \subseteq \{0,1\}^n$ and a bilinear loss function $\ell(\mathbf{x}, \mathbf{c}) = \mathbf{c}^\top \mathbf{x}$. Also, for simplicity the support set (\ref{eq: support constraints}) is assumed to be component-wise interval and given by: 
\begin{equation} \label{eq: support constraints experiments} 
	\mathcal{S}_0 := \mathcal{S}^{(bf)}_0 = \Big\{\mathbf{c}' \in \mathbb{R}^{n}: \mathbf{0} \leq \mathbf{c}' \leq \mathbf{1}\Big\}.
\end{equation}	
 For all problems, except for the maximum coverage problem, we make the assumption that any decision $\mathbf{x} \in X$ satisfies $\Vert \mathbf{x} \Vert_1 = h$ for some $h \in \mathbb{Z}_{++}$. 

\textit{\textbf{Sorting problem (SP).}} In the sorting problem there are $n$ items. With each item $a \in \mathcal{A} = \{1, \ldots, n\}$ we associate a nonnegative cost $c_a$ and attempt to minimize the total expected cost of $h \in \mathbb{Z}_{++}$ selected items. Formally, we set 
\begin{equation} \label{eq: feasible set for sorting problem} \nonumber
X^{(SP)} = \Big\{\mathbf{x} \in \{0, 1\}^{n}: \Vert \mathbf{x} \Vert_1 = h \Big\}.
\end{equation}

In addition to the support constraints (\ref{eq: support constraints experiments}), we introduce component-wise interval linear data constraints (\ref{eq: data uncertainty}) given by
\begin{equation}
\hat{\mathbf{c}}^{\mbox{\tiny (\itshape k\upshape)}} \in \mathcal{S}^{(SP)}_k := \Big\{\mathbf{c}' \in \mathbb{R}^{n}: \max\{\hat{c}^{\mbox{\tiny (\itshape k\upshape)}}_a - \delta^{\mbox{\tiny (\itshape k\upshape)}}_a, 0\} \leq c'_a \leq \min\{\hat{c}^{\mbox{\tiny (\itshape k\upshape)}}_a + \delta^{\mbox{\tiny (\itshape k\upshape)}}_a, 1\}\ \; \forall a \in \mathcal{A} \Big\}, \quad \forall k \in \mathcal{K}. \label{eq: data constraints sorting problem}
\end{equation} 
In other words, the random observations $\hat{c}^{\mbox{\tiny (\itshape k\upshape)}}_a$ for each $a \in \mathcal{A}$ are subject to noise, whose magnitude is controlled by additional parameters $\delta^{\mbox{\tiny (\itshape k\upshape)}}_a \in [0, 1]$. 

We recall that according to Theorem \ref{theorem 2} the resulting three-level problem (\ref{three level formulation}) with interval uncertainty can be solved efficiently since the sorting problem is polynomially solvable. As a remark, we note that the stochastic sorting problem with $h = 1$ can be viewed as a static version of the well-known online multi-armed bandit problem \cite{Auer2002}. 

\textbf{\textit{Shortest path problem (SPP).}} As a second class of problems, we consider the SPP in a fully-connected acyclic layered graph with $ h - 1 \in \mathbb{Z}_{++}$ intermediate layers and $r \in \mathbb{Z}_{++}$ nodes at each layer. The first and the last layers consist of unique nodes, which are the source and the destination nodes, respectively, and each path contains exactly $h$ arcs; see, e.g., Figure \ref{fig: layered graph}. In particular, the indices $\mathcal{A} = \{1, \ldots, n\}$ are related to the set of directed arcs with random costs $c_a$, $a \in \mathcal{A}$, whereas a decision $\mathbf{x} \in \{0, 1\}^n$ encodes a simple path between the source and the destination nodes. Hence, $X^{(SPP)}$ is given by standard path flow constraints \cite{Ahuja1988}, whose explicit form is omitted for brevity. %and the decision-maker's objective is to select a path with the least possible expected cost.

\begin{figure}
	\centering
	\begin{subfigure}{0.49\textwidth} 	\centering
		\begin{tikzpicture}[scale=0.6,transform shape]
			\Vertex[x=0,y=0]{1}
			\Vertex[x=3,y=3]{2}
			\Vertex[x=3,y=0]{3}
			\Vertex[x=3,y=-3]{4}
			\Vertex[x=6,y=3]{5}
			\Vertex[x=6,y=0]{6}
			\Vertex[x=6,y=-3]{7}
			\Vertex[x=9,y=0]{8}
			\tikzstyle{LabelStyle}=[fill=white,sloped]
			\tikzstyle{EdgeStyle}=[post]
			\Edge(1)(2)
			\Edge(2)(6)
			\Edge(7)(8)
			\Edge(1)(3)
			\Edge(1)(4)
			\Edge(5)(8)
			\Edge(6)(8)
			\Edge(2)(5)
			\Edge(2)(7)
			\Edge(3)(5)
			\Edge(3)(6)
			\Edge(3)(7)
			\Edge(4)(5)
			\Edge(4)(6)
			\Edge(4)(7)
		\end{tikzpicture}
		\caption{\footnotesize \centering A fully-connected layered graph with $h = 3$ and $r = 3$ for the SPP. }
		\label{fig: layered graph}
	\end{subfigure}
	\hfill
	\begin{subfigure}{0.49\textwidth} 	\centering
			\begin{tikzpicture}[scale=0.6,transform shape]
			\Vertex[x=-3,y=3,L =$\mathcal{A}_1$]{1 }
			\Vertex[x=0,y=3,L =$\mathcal{A}_2$]{2 }
			\Vertex[x=3,y=3,L =$\mathcal{A}_3$]{3 }
			\Vertex[x=-3,y=-3]{1}
			\Vertex[x=-1,y=-3]{2}
			\Vertex[x=1,y=-3]{3}
			\Vertex[x=3,y=-3]{4}
			\tikzstyle{LabelStyle}=[fill=white,sloped]
			\tikzstyle{EdgeStyle}=[post]
			\Edge(1 )(1)
			\Edge(1 )(2)
			\Edge(2 )(1)
			\Edge(2 )(3)
			\Edge(3 )(2)
			\Edge(3 )(3)
			\Edge(3 )(4)
		\end{tikzpicture}
		\caption{\footnotesize \centering An instance of the MCP with $n_1 = 3$, $n_2 = 4$, $\mathcal{A}_1 = \{1, 2\}$, $ \mathcal{A}_2 = \{1, 3\}$ and $\mathcal{A}_3 = \{2, 3, 4\}$. }
		\label{fig: bipartite graph}
	\end{subfigure}
	\label{fig: graphs}
	
\end{figure}

As briefly outlined in Section \ref{subsec: our approach and contributions}, we may consider two problem-specific forms of data uncertainty, semi-bandit and bandit feedback. Thus, for a set of $K$ directed paths, say $\mathcal{P} = \{P^{\mbox{\tiny (\itshape k\upshape)}}, k \in \mathcal{K}\}$, from the source to the destination node, we assume that either (\textit{i}) the decision-maker observes the cost of each arc $a \in P^{\mbox{\tiny (\itshape k\upshape)}}$ for $k \in \mathcal{K}$ or (\textit{ii}) observes only the total cost of each path $P^{\mbox{\tiny (\itshape k\upshape)}} \in \mathcal{P}$. Formally, for each $k \in \mathcal{K}$ we introduce the scenarios of semi-bandit and bandit feedback, respectively, as follows: 
\begin{subequations} \label{eq: data uncertainty spp}
\begin{align}
&\mathcal{S}^{(SPP_1)}_k := \Big\{\mathbf{c}' \in \mathbb{R}^{n}: c'_a = c^{\mbox{\tiny (\itshape k\upshape)}}_a \quad \forall a \in P^{\mbox{\tiny (\itshape k\upshape)}} \Big\} \cap \mathcal{S}_0, \label{eq: data uncertainty sbf} \\
& \mathcal{S}^{(SPP_2)}_k := \Big\{\mathbf{c}' \in \mathbb{R}^{n}: \sum_{a \in P^{\mbox{\tiny (\itshape k\upshape)}}} c'_a = S^{\mbox{\tiny (\itshape k\upshape)}} \Big\} \cap \mathcal{S}_0, \label{eq: data uncertainty bf}
\end{align}
\end{subequations}
where $\mathcal{S}_0$ is defined by (\ref{eq: support constraints experiments}); $c^{\mbox{\tiny (\itshape k\upshape)}}_a \in \mathbb{R}_{+}$ and $S^{\mbox{\tiny (\itshape k\upshape)}} \in \mathbb{R}_{+}$ are the costs of arcs contained in $P^{\mbox{\tiny (\itshape k\upshape)}}$ and the total cost of $P^{\mbox{\tiny (\itshape k\upshape)}}$, respectively. A way to select the set of paths $\mathcal{P}$ is described later in Section \ref{subsec: results and discussion}

\textit{\textbf{Maximum coverage problem (MCP).}} In the MCP we a given an integer number $\tilde{h} \in \mathbb{Z}_{++}$ and a collection of subsets $\mathcal{A}_i \subseteq \mathcal{A} = \{1, \ldots, n_1\}$, $i \in \{1, \ldots, n_2\}$. The goal is to find a subcollection, whose cardinality does not exceed $\tilde{h}$ and the expected cost of covered items in $\mathcal{A}$ is maximized. Formally, we introduce binary variables $x_a \in \{0, 1\}$ and $y_i \in \{0, 1\}$, which indicate, respectively, whether the item $a \in \mathcal{A}$ is covered or not and whether the subset $\mathcal{A}_i \subseteq \mathcal{A}$, $i \in \{1, \ldots, n_2\}$, is selected or not. Then, the set of feasible decisions can be defined as:
\begin{equation} \nonumber
X^{(MCP)} = \Big\{(\mathbf{x}, \mathbf{y}) \in \{0, 1\}^{n_1 + n_2}: \sum_{i = 1}^{n_2} y_i \leq \tilde{h}, \; \sum_{i : \; a \in \mathcal{A}_i} y_i \geq x_a \quad \forall a \in \mathcal{A} \Big\}.	
\end{equation}
Furthermore, the loss function is given by $\ell(\mathbf{x}, \mathbf{c}) = -\mathbf{c}^\top \mathbf{x}$, where $c_a$ is a random cost of $a \in \mathcal{A}$. In contrast to the sorting and the shortest path problems, the MCP is known to be $NP$-hard \cite{Garey2002} and the assumption that each decision $(\mathbf{x}, \mathbf{y}) \in X$ has a fixed number of non-zero elements is not satisfied. Any instance of the MCP can be represented using a bipartite graph as illustrated in Figure \ref{fig: bipartite graph}. 

Next, similar to the SPP, we may define a set of decisions $\widetilde{\mathcal{P}} = \{\tilde{P}^{\mbox{\tiny (\itshape k\upshape)}}, k \in \mathcal{K}\}$, where $\tilde{P}^{\mbox{\tiny (\itshape k\upshape)}}$ is a subcollection of $\tilde{h}$ subsets selected from $\mathcal{A}_1, \ldots, \mathcal{A}_{n_2}$. For a given $\tilde{P}^{\mbox{\tiny (\itshape k\upshape)}}$ we can identify the covered elements in $\mathcal{A}$, i.e., for each $a \in \mathcal{A}$
\begin{equation} \nonumber
	x^{\mbox{\tiny (\itshape k\upshape)}}_a = \begin{cases}
		1, \mbox{ if } a \mbox{ is covered by } \tilde{P}^{\mbox{\tiny (\itshape k\upshape)}}, \\
		0, \mbox{ otherwise.}
	\end{cases}
\end{equation}
Then, for each subcollection of sets $\tilde{P}^{\mbox{\tiny \upshape (\itshape k\upshape)}}$, $k \in \mathcal{K}$, we define the scenarios of semi-bandit and bandit feedback, respectively, as follows: 
\begin{subequations} \label{eq: data uncertainty mcp}
	\begin{align} 
		&\mathcal{S}^{(MCP_1)}_k := \Big\{\mathbf{c}' \in \mathbb{R}^{n}: c'_a = \tilde{c}^{\mbox{\tiny (\itshape k\upshape)}}_a \quad \forall a: \; x^{\mbox{\tiny (\itshape k\upshape)}}_a = 1 \Big\} \cap \mathcal{S}_0, \label{eq: data uncertainty mcp sbf} \\
		& \mathcal{S}^{(MCP_2)}_k := \Big\{\mathbf{c}' \in \mathbb{R}^{n}: \sum_{a: \; x^{\mbox{\tiny (\itshape k\upshape)}}_a = 1} c'_a = \tilde{S}^{\mbox{\tiny (\itshape k\upshape)}} \Big\} \cap \mathcal{S}_0. \label{eq: data uncertainty mcp bf}
	\end{align}
\end{subequations}
Here, $\tilde{c}^{\mbox{\tiny (\itshape k\upshape)}}_a \in \mathbb{R}_{+}$ and $\tilde{S}^{\mbox{\tiny (\itshape k\upshape)}} \in \mathbb{R}_+$ are the costs of items covered by $\tilde{P}^{\mbox{\tiny (\itshape k\upshape)}}$ and their total cost, respectively. 

The considered forms of data uncertainty both for the SPP and the MCP are motivated by the related online problem settings described in \cite{Chen2013}. However, the online version of MCP studied in \cite{Chen2013} is essentially an unweighted MCP with uncertain "activation probabilities" for subset-item pairs and a non-linear loss function. Taking into account Assumption \textbf{A2}, we present a modified version of MCP that maintains the linearity of the loss function $\ell(\mathbf{x}, \mathbf{c})$.

\textbf{Nominal distribution.} Taking into account the specific form of the support set (\ref{eq: support constraints experiments}), there is no initial information regarding the dependence between the components of the cost vector $\mathbf{c}$. Therefore, to simplify the analysis, we consider joint distributions $\mathbb{Q}^* \in \mathcal{Q}_0(\mathcal{S}_0)$ with independent components. More precisely, we assume that the costs $c_a$ for each $a \in \mathcal{A}$ are governed by a standard beta distribution with parameters $\alpha_a, \beta_a \in \mathbb{R}_{> 0}$ and a support given by $[0, 1]$. 
The parameters $\alpha_a$ and $\beta_a$ can be defined using the mean, $m_a$, and the standard deviation, $\sigma_a$, of $c_a$, i.e.,
\begin{equation} \label{eq: parameters of beta distribution}
\begin{gathered}
\alpha_a = \frac{m_a^2(1 - m_a)}{\sigma^2_a} - m_a, \quad
\beta_a = \alpha_a(\frac{1}{m_a} - 1),
\end{gathered}
\end{equation}
see, e.g., \cite{Gupta2004}. In all experiments we set $\sigma_{a} = 0.125$, $a \in \mathcal{A}$, and select $m_{a}$ uniformly at random from the interval
$$\Big(\frac{1}{2}(1 - \sqrt{1 - 4\sigma_{a}^2}), \frac{1}{2}(1 + \sqrt{1 - 4\sigma_{a}^2})\Big).$$
The latter condition guarantees that a beta distribution defined by (\ref{eq: parameters of beta distribution}) exists, i.e., $\alpha_{a}, \beta_{a} > 0$. Finally, the joint distribution $\mathbb{Q}^*$ is defined as a product of the associated marginal distributions. 

\textbf{Wasserstein radius.} %Several existing approaches to estimating the Wasserstein radius $\varepsilon_K$ for DRO problems with \textit{complete data} are discussed in \cite{Esfahani2018, Kuhn2019}. 
In fact, by selecting an appropriate value of the Wasserstein radius $\varepsilon_K$, we would like to guarantee that the nominal distribution $\mathbb{Q}^*$ belongs to the Wasserstein ball centered at the empirical distribution of the data or, equivalently, $W^1(\widehat{\mathbb{Q}}_K, \mathbb{Q}^*) \leq \varepsilon_K$, with high probability. Although the data set $\widehat{\mathbf{C}}$ in our problem setting is incomplete, one may argue that for any fixed $\mathbf{x} \in X$ %In particular, if the data set $\widehat{\mathbf{C}}$ is fixed, then it is sufficient to guarantee that for any decision $\mathbf{x} \in X$ and a prescribed confidence level $\eta \in (0, 1)$ 
\begin{equation} \label{eq: finite sample guarantee}
\begin{gathered}
\Pr\Big\{\mathbb{E}_{\mathbb{Q}^*}\{\ell(\mathbf{x}, \mathbf{c})\} > \max_{ \widehat{\mathbf{C}} } \big\{ \max_{\mathbb{Q} \in \mathcal{Q}(\widehat{\mathbf{C}})} \mathbb{E}_\mathbb{Q} \{\ell(\mathbf{x}, \mathbf{c})\}: \hat{\mathbf{c}}^{\mbox{\tiny (\itshape k\upshape)}} \in \mathcal{S}_k \; \forall k \in \mathcal{K} \big\}\Big\} \leq \\ 
\Pr\Big\{\mathbb{E}_{\mathbb{Q}^*}\{\ell(\mathbf{x}, \mathbf{c})\} > \max_{\mathbb{Q} \in \mathcal{Q}(\widetilde{\mathbf{C}})} \mathbb{E}_\mathbb{Q} \{\ell(\mathbf{x}, \mathbf{c})\}\Big\} \leq \Pr\Big\{W^1(\widetilde{\mathbb{Q}}_K, \mathbb{Q}^*) > \varepsilon_K \Big\},
\end{gathered}
\end{equation} 
where $\widetilde{\mathbf{C}}$ and $\widetilde{\mathbb{Q}}_K$ are the nominal data set obtained from $\mathbb{Q}^*$ and its empirical distribution, respectively. Thus, if we guarantee that
$$\Pr\Big\{W^1(\widetilde{\mathbb{Q}}_K, \mathbb{Q}^*) > \varepsilon_K \Big\} \leq \eta$$ for a sufficiently small $\eta \in (0, 1) $, then equation (\ref{eq: finite sample guarantee}) asserts that with probability of at least $1 - \eta$ the worst-case expected loss in our setting provides an \textit{upper bound} on
the nominal expected loss. %; see, e.g., \cite{Esfahani2018, Kuhn2019} for more~detail. 

We demonstrate that in the case of interval support constraints (\ref{eq: support constraints experiments}) used in our experiments, the Wasserstein radius can be defined as $\varepsilon_K = \frac{\gamma}{\sqrt{K}}$ for some $\gamma \in \mathbb{R}_{++}$. Indeed, by using Lemma \ref{lemma 1} with $\widehat{\mathbf{C}} = \widetilde{\mathbf{C}}$ and Hoeffding inequality \cite{Hoeffding1963}, we observe that: 
\begin{align} \allowdisplaybreaks \label{eq: Hoeffding}
\Pr\Big\{\mathbb{E}_{\mathbb{Q}^*}\{\ell(\mathbf{x}, \mathbf{c})\} > \max_{\mathbb{Q} \in \mathcal{Q}(\widetilde{\mathbf{C}})} \mathbb{E}_\mathbb{Q} \{\ell(\mathbf{x}, \mathbf{c})\}\Big\} & = \Pr\Big\{\mathbb{E}_{\mathbb{Q}^*}\{\mathbf{c}^\top \mathbf{x}\} > \min\big\{\frac{1}{K}\sum_{k = 1}^K \tilde{\mathbf{c}}^{\mbox{\tiny \upshape (\itshape k\upshape)}\top} \mathbf{x} + \varepsilon_K; \mathbf{1}^\top \mathbf{x} \big\}\Big\} \leq \nonumber \\
\Pr\Big\{\mathbb{E}_{\mathbb{Q}^*}\{\mathbf{c}^\top \mathbf{x}\} > \frac{1}{K}\sum_{k = 1}^K \tilde{\mathbf{c}}^{\mbox{\tiny \upshape (\itshape k\upshape)}\top} \mathbf{x} + & \varepsilon_K \Big\} + \Pr\Big\{\mathbb{E}_{\mathbb{Q}^*}\{\mathbf{c}^\top \mathbf{x}\} > \mathbf{1}^\top \mathbf{x} \Big\} = \\
\Pr\Big\{\mathbb{E}_{\mathbb{Q}^*}\{\mathbf{c}^\top \mathbf{x}\} > \frac{1}{K}\sum_{k = 1}^K \tilde{\mathbf{c}}^{\mbox{\tiny \upshape (\itshape k\upshape)}\top} \mathbf{x} + & \varepsilon_K \Big\} \leq \exp\Big(- \frac{2 K \varepsilon^2_K}{\big( \max_{\mathbf{x} \in X} \Vert \mathbf{x} \Vert_1 \big) ^2} \Big). \nonumber
\end{align}
In particular, the equality follows from the fact that $\mathbb{E}_{\mathbb{Q}^*}\{\mathbf{c}^\top \mathbf{x}\} > \mathbf{1}^\top \mathbf{x}$ with zero probability; recall~(\ref{eq: support constraints experiments}). 
Consequently, by setting $\varepsilon_K$ in the order of $K^{-\frac{1}{2}}$ we may guarantee that the inequality 
$$\mathbb{E}_{\mathbb{Q}^*}\{\ell(\mathbf{x}, \mathbf{c})\} \leq \max_{ \widehat{\mathbf{C}} } \big\{ \max_{\mathbb{Q} \in \mathcal{Q}(\widehat{\mathbf{C}})} \mathbb{E}_\mathbb{Q} \{\ell(\mathbf{x}, \mathbf{c})\}: \hat{\mathbf{c}}^{\mbox{\tiny (\itshape k\upshape)}} \in \mathcal{S}_k \; \forall k \in \mathcal{K} \big\}$$
holds with high probability. 

Unfortunately, the choice of $\varepsilon_K$ becomes substantially more complicated when dealing with a general polyhedral support set $\mathcal{S}_0$. First, we note that the standard assumption that the nominal distribution $\mathbb{Q}^*$ is light-tailed \cite{Esfahani2018} is readily satisfied in our setting as $\mathcal{S}_0$ is bounded. 
To the best of our knowledge, the most effective finite-sample guarantees for the Wasserstein distance are developed by Gao \cite{Gao2022b}. In particular, it is shown that, under some additional assumptions about the loss function $\ell(\mathbf{x}, \mathbf{c})$, $\varepsilon_K$ can be selected in the order of $K^{-\frac{1}{2}}$ up to a logarithmic factor (see Corollary 4 in~\cite{Gao2022b} for a particular case of type-1 Wasserstein balls).
At the same time, the obtained estimates of $\varepsilon_K$ implicitly depend on parameters of $\mathbb{Q}^*$ and, therefore, their practical implementation is quite limited. We leave a further discussion of this issue as a possible direction for future research. 

\textbf{Computational settings.} All experiments are performed on a PC with CPU i5-7200U and RAM 8 GB. The linear MIP reformulation (\ref{mixed-integer programming reformulation}) as well as deterministic versions of the considered combinatorial optimization problems are solved in Java with CPLEX 20.1.

\subsection{Results and discussion} \label{subsec: results and discussion}
\textbf{Experiments for the sorting problem.}
In the first set of experiments we consider the sorting problem with interval uncertainty defined by (\ref{eq: data constraints sorting problem}) and analyze the quality of distributionally robust solutions as a function of various parameters of the problem. In view of Theorem \ref{theorem 2}, the three-level problem (\ref{three level formulation}) with interval uncertainty can be resolved via finding the optimal objective function values of two deterministic sorting problems, (\ref{mixed-integer programming reformulation 1 interval case}) and (\ref{mixed-integer programming reformulation 2 interval case}). In all experiments we set the number of items $n = 50$ and compute the average relative loss (\ref{eq: nominal relative loss}) with mean absolute deviations (MADs) for $100$ randomly generated test instances. 

First, let $K = 50$ and assume that that each element $c^{\mbox{\tiny (\itshape k\upshape)}}_a$ of the data set (\ref{eq: data set}) is subject to interval uncertainty with a probability $p_a \in [0,1]$, where the parameters $p_a$, $a \in \mathcal{A}$, are selected uniformly at random from the interval $[0,1]$. More specifically, for each $a \in \mathcal{A}$ and $k \in \mathcal{K}$ we set
\begin{equation} \nonumber
\delta^{\mbox{\tiny (\itshape k\upshape)}}_a = \begin{cases} \delta \mbox{ with probability } p_a,\\
0 \mbox{ with probability } 1 - p_a. \end{cases}
\end{equation}
We consider the nominal relative loss (\ref{eq: nominal relative loss}) as a function of the noise level $\delta$ (Figure \ref{fig: dependence on delta}) and the number of selected items $h$ (Figure \ref{fig: dependence on h}). Finally, by varying $\gamma$ for fixed $\delta$ and $h$ we explore how the nominal relative loss scales in the Wasserstein radius $\varepsilon_K = \frac{\gamma}{\sqrt{K}}$; see Figure \ref{fig: dependence on epsilon}. 

We make the following observations:
\begin{itemize}
	\item From Figure~\ref{fig: dependence on delta} we observe that the nominal relative loss, $\rho$, tends to increase with the increase of $\delta$. This fact is rather intuitive, as with the increase of $\delta$ we need to resort to more conservative decisions resulting in poorer out-of-sample performance. 
\begin{figure}
	 \centering
\begin{subfigure}{0.3\textwidth}		
\begin{tikzpicture}[scale = 0.6]
\begin{axis}[
	xlabel=\footnotesize {$\qquad \quad \delta$},
ylabel={ \footnotesize $\qquad$ Nominal relative loss},
xmin=0, xmax=1,
ymin=0, ymax=3.5,
xtick={0.2,0.4,0.6,0.8,1},
ytick={0,1,1.5,2,2.5,3,3.5},
legend pos=north west,
ymajorgrids=true,
grid style=dashed,
]

\addplot[name path=f1,
color=blue,
mark=square,
]
coordinates {
	(0.00, 1.03) 
	(0.10, 1.11) 
	(0.20, 1.28) 
	(0.30, 1.48) 
	(0.40, 1.73) 
	(0.50, 1.88) 
	(0.60, 2.02) 
	(0.70, 2.15) 
	(0.80, 2.34) 
	(0.90, 2.47)
};
\addplot[name path=f2,
color=blue,
style=dashed,
mark=-,
]
coordinates {
	(0.00, 0.99) 
	(0.10, 1.02) 
	(0.20, 1.1) 
	(0.30, 1.22) 
	(0.40, 1.33) 
	(0.50, 1.41) 
	(0.60, 1.54) 
	(0.70, 1.59) 
	(0.80, 1.68) 
	(0.90, 1.75)
};
\addplot[name path=f3,
color=blue,
style=dashed,
mark=-,
]
coordinates {
	(0.00, 1.07) 
	(0.10, 1.2) 
	(0.20, 1.46) 
	(0.30, 1.74) 
	(0.40, 2.13) 
	(0.50, 2.35) 
	(0.60, 2.5) 
	(0.70, 2.71) 
	(0.80, 3) 
	(0.90, 3.19)
};

\addplot [
thick,
color=blue,
fill=blue, 
fill opacity=0.05
]
fill between[
of=f3 and f2,
soft clip={domain=0:1},
];
\end{axis}
\end{tikzpicture}
\caption{\footnotesize $K = 50$, $h = 5$ and $\varepsilon_K = 1$.}
\label{fig: dependence on delta}
\end{subfigure}
\hfill
\begin{subfigure}{0.3\textwidth}
		\begin{tikzpicture}[scale = 0.6]
		\begin{axis}[
		xlabel=\footnotesize {$\qquad \quad h$},
	ylabel={ \footnotesize $\qquad$ Nominal relative loss},
		xmin=0, xmax=10,
		ymin=0, ymax=16.5,
		xtick={2,4,6,8,10},
		ytick={1,5,10,15},
		legend pos=north west,
		ymajorgrids=true,
		grid style=dashed,
		]
		
		\addplot[name path=f1,
		color=blue,
		mark=square,
		]
		coordinates {
			(1,16.13)
			(2,1.66)
			(3,1.45)
			(4,1.33)
			(5,1.26)
			(6,1.22)
			(7,1.16)
			(8,1.12)
			(9,1.09)
			(10,1.07)
		};
		\addplot[name path=f2,
		color=blue,
		style=dashed,
		mark=-,
		]
		coordinates {
			(1,7.34)
			(2,1.09)
			(3,1.08)
			(4,1.09)
			(5,1.07)
			(6,1.08)
			(7,1.06)
			(8,1.04)
			(9,1.04)
			(10,1.02)
		};
		\addplot[name path=f3,
		color=blue,
		style=dashed,
		mark=-,
		]
		coordinates {
			(1,24.93)
			(2,2.24)
			(3,1.83)
			(4,1.58)
			(5,1.44)
			(6,1.36)
			(7,1.26)
			(8,1.20)
			(9,1.14)
			(10,1.12)
		};
		
		\addplot [
		thick,
		color=blue,
		fill=blue, 
		fill opacity=0.05
		]
		fill between[
		of=f3 and f2,
		soft clip={domain=0:1},
		];
		\end{axis}
		\end{tikzpicture}
		\caption{\footnotesize \centering $K = 50$, $\delta = 0.2$ and $\varepsilon_K = 1$.}
		\label{fig: dependence on h}
\end{subfigure}
\hfill 
\begin{subfigure}{0.3\textwidth}
	\centering
	\begin{tikzpicture}[scale = 0.6]
	\begin{axis}[
		xlabel=\footnotesize {$\qquad \quad \gamma$},
	ylabel={ \footnotesize $\qquad$ Nominal relative loss},
	xmin=0, xmax=55,
	ymin=0, ymax=12,
	xtick={10,20,30,40,50},
	ytick={1,5,10},
	legend style={at={(0.45,1.5)},anchor=north},
	ymajorgrids=true,
	grid style=dashed,
	]
	
	\addplot[name path=f1,
	color=blue,
	mark=square,
	]
	coordinates {
		(5.00,	1.04)
		(10.00,	1.04)
		(15.00,	1.04)
		(20.00,	1.04)
		(25.00,	1.04)
		(30.00,	1.04)
		(35.00,	8.54)
		(40.00,	8.54)
		(45.00,	8.54)
		(50.00,	8.54)
	};
	
	\addplot[name path=f4,
	color=red,
	mark=triangle ,
	]
	coordinates {
		(5.00,	1.33)
		(10.00,	1.33)
		(15.00,	1.33)
		(20.00,	1.33)
		(25.00,	1.33)
		(30.00,	3.63)
		(35.00,	8.54)
		(40.00,	8.54)
		(45.00,	8.54)
		(50.00,	8.54)
	};
	
	\addplot[name path=f2,
	color=blue,
	style=dashed,
	mark=-,
	]
	coordinates {
		(5.00,	0.99)
		(10.00,	0.99)
		(15.00,	0.99)
		(20.00,	0.99)
		(25.00,	0.99)
		(30.00,	0.99)
		(35.00,	5.81)
		(40.00,	5.81)
		(45.00,	5.81)
		(50.00,	5.81)
	};
	\addplot[name path=f3,
	color=blue,
	style=dashed,
	mark=-,
	]
	coordinates {
		(5.00,	1.08)
		(10.00,	1.08)
		(15.00,	1.08)
		(20.00,	1.08)
		(25.00,	1.08)
		(30.00,	1.08)
		(35.00,	11.26)
		(40.00,	11.26)
		(45.00,	11.26)
		(50.00,	11.26)
	};
	
	\addplot[name path=f5,
	color=red,
	style=dashed,
	mark=-,
	]
	coordinates {
		(5.00,	1.09)
		(10.00,	1.09)
		(15.00,	1.09)
		(20.00,	1.09)
		(25.00,	1.09)
		(30.00,	0.77)
		(35.00,	5.81)
		(40.00,	5.81)
		(45.00,	5.81)
		(50.00,	5.81)
	};
	\addplot[name path=f6,
	color=red,
	style=dashed,
	mark=-,
	]
	coordinates {
		(5.00,	1.56)
		(10.00,	1.56)
		(15.00,	1.56)
		(20.00,	1.56)
		(25.00,	1.56)
		(30.00,	6.48)
		(35.00,	11.26)
		(40.00,	11.26)
		(45.00,	11.26)
		(50.00,	11.26)
	};
	
	\addplot [
	thick,
	color=blue,
	fill=blue, 
	fill opacity=0.05
	]
	fill between[
	of=f3 and f2,
	soft clip={domain=5:50},
	];
	
	\addplot [
	thick,
	color=red,
	fill=red, 
	fill opacity=0.05
	]
	fill between[
	of=f6 and f5,
	soft clip={domain=5:50},
	];
	\legend{\footnotesize $\delta=0$, \footnotesize $\delta=0.2$}
	\end{axis}
	\end{tikzpicture}
\caption{\footnotesize $K = 50$, $h = 5$ and $\delta = 0$ (complete data) or $\delta = 0.2$ (noisy data).}
\label{fig: dependence on epsilon}
\end{subfigure} 
\caption{\footnotesize We report the average relative loss (\ref{eq: nominal relative loss}) with MADs as a function of $\delta$ (a), $h$ (b) and $\gamma$ (c) for $100$ random test instances.}
\end{figure}
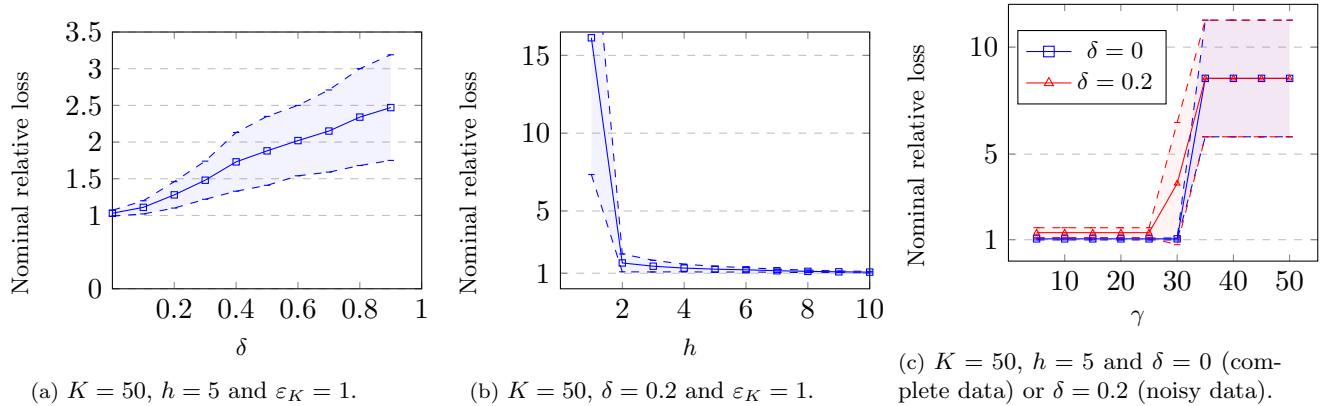 

	\item From Figure~\ref{fig: dependence on h} we observe that the average value of $\rho$ decreases in $h$. This observation can be justified by Theorem \ref{theorem 2} and the fact that $\varepsilon_K = 1$. Specifically, when $h=1$, the optimal objective function value of (\ref{mixed-integer programming reformulation 1 interval case}) is always at least 1, while the optimal objective function value of (\ref{mixed-integer programming reformulation 2 interval case}) equals 1. Hence, for $h = 1$ we always resolve ties to the most conservative solution. The larger $h$, the more often an optimal solution of (\ref{three level formulation}) is provided by (\ref{mixed-integer programming reformulation 1 interval case}), which implies the~result. 
	\item The behavior of $\rho$ exhibits two distinct regimes depending on the value of $\gamma$ or $\varepsilon_K$, as depicted in Figure \ref{fig: dependence on epsilon}. This observation can also be explained by Theorem \ref{theorem 2}, i.e., an optimal solution of (\ref{three level formulation}) alternates between a solution obtained via robust sample average approximation and the worst-case solution corresponding to the vector $\mathbf{u} = \mathbf{1}$ of upper bounds. 
	\item From Figure \ref{fig: dependence on epsilon} we also conclude that the average values of $\rho$ for $\delta = 0$ are smaller than those for $\delta = 0.2$; this fact follows from the first observation (Figure \ref{fig: dependence on delta}). 	 
\end{itemize}

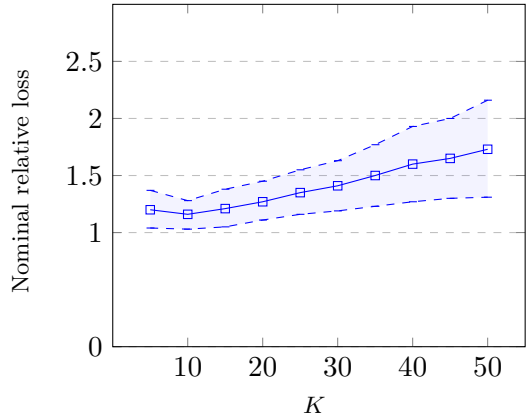
\begin{figure}	
	\centering
	\begin{tikzpicture}[scale = 0.8]
	\begin{axis}[
	xlabel=\footnotesize {$\qquad \quad K$},
	ylabel={ \footnotesize $\qquad$ Nominal relative loss},
	xmin=0, xmax=55,
	ymin=0, ymax=3,
	xtick={10,20,30,40,50},
	ytick={0,1,1.5,2,2.5},
	legend pos=north west,
	ymajorgrids=true,
	grid style=dashed,
	]
	
	\addplot[name path=f1,
	color=blue,
	mark=square,
	]
	coordinates {
		(5,	1.20)
		(10,	1.16)
		(15,	1.21)
		(20,	1.27)
		(25,	1.35)
		(30,	1.41)
		(35,	1.50)
		(40,	1.60)
		(45,	1.65)
		(50,	1.73)
	};
	
	\addplot[name path=f2,
	color=blue,
	style=dashed,
	mark=-,
	]
	coordinates {
		(5,	1.04)
		(10,	1.03)
		(15,	1.05)
		(20,	1.11)
		(25,	1.16)
		(30,	1.19)
		(35,	1.23)
		(40,	1.27)
		(45,	1.30)
		(50,	1.31)
		};
	\addplot[name path=f3,
	color=blue,
	style=dashed,
	mark=-,
	]
	coordinates {
		(5,	1.37)
		(10,	1.28)
		(15,	1.38)
		(20,	1.45)
		(25,	1.55)
		(30,	1.63)
		(35,	1.77)
		(40,	1.93)
		(45,	2.00)
		(50,	2.16)
	};
	
		\addplot [
	thick,
	color=blue,
	fill=blue, 
	fill opacity=0.05
	]
	fill between[
	of=f3 and f2,
	soft clip={domain=5:50},
	];
 \end{axis}
	\end{tikzpicture}
	\caption{\footnotesize Let $K_{max} = 50$, $\gamma = \sqrt{K_{max}}$ and $h = 5$. We report the average relative loss (\ref{eq: nominal relative loss}) with MADs as a function of $K$ for $100$ random test instances.}
	\label{fig: dependence on sample size}
\end{figure}

In the next experiment we consider the sorting problem in a similar setting but vary the level of noise and the sample size. More specifically, we assume that for each $a \in \mathcal{A}$ and $k \in \mathcal{K}$
\begin{equation} \nonumber
\delta^{\mbox{\tiny (\itshape k\upshape)}}_a = \begin{cases} \frac{k - 1}{K_{max}} \mbox{ with probability } p_a,\\
0 \mbox{ with probability } 1 - p_a, \end{cases}
\end{equation}
where $K_{max} = 50$ is the maximal considered sample size. Formally, each component of the cost vector~$\mathbf{c}$ is also perturbed with probability $p_a \in [0, 1]$, $a \in \mathcal{A}$, and the magnitude of noise increases with the increase of $k \in \mathcal{K}$. In contrast to the previous experiments, we change the sample size $K$ and, in particular, set $\gamma = \sqrt{K_{max}}$ and $\varepsilon_K = \sqrt{\frac{K_{max}}{K}}$. The results are reported in Figure \ref{fig: dependence on sample size}. 

\looseness-1 We observe that the nominal relative loss (\ref{eq: nominal relative loss}) first tends to decrease and then increases as a function of $K$. The intuition behind this trend can be explained as follows. On the one hand, as $K$ increases, the Wasserstein radius $\varepsilon_K$ decreases, and by the law of large numbers, we can better approximate the expected loss in (\ref{stochastic programming problem}) by its sample mean; recall Theorem \ref{theorem 2}. On the other hand, as $K$ increases, we subsequently append random samples with a higher magnitude of noise that clearly leads to a misspecification of the optimal solution of (\ref{stochastic programming problem}) for sufficiently large values of $K$. We conclude that using incomplete data with a reasonable magnitude of noise (not exceeding $20\%$ of its range) may improve the model's out-of-sample performance. 

\begin{algorithm}
	\DontPrintSemicolon
	\footnotesize For each arc $a \in \mathcal{A}$ maintain: (1) variable $T_a$ as the
	total number of times the arc $a$ is observed so far; (2)
	variable $\hat{\mu}_a$ as the empirical mean of the random cost $c_a$ (initially we set $T_a = 0$ and $\hat{\mu}_a = 0$ for each $a \in \mathcal{A}$); \;
	\textbf{for} $k \in \{1, \ldots, K\}$: \;
	\Begin{
		find the shortest path $P^{\mbox{\tiny (\itshape k\upshape)}}$ with respect to the costs $\hat{c}_a = \max\{\hat{\mu}_a - \sqrt{\frac{3 \ln k}{2T_a}}, 0\}$, $a \in \mathcal{A}$;\;
		update $T_a$ and $\hat{\mu}_a$ for each $a \in P^{\mbox{\tiny (\itshape k\upshape)}}$;
	}
	\textbf{return} $\mathcal{P} = \{P^{\mbox{\tiny (\itshape k\upshape)}}, k \in \mathcal{K}\}$. 
	\caption{CUCB algorithm for generating the set of paths \cite{Wang2017}.}
	\label{alg: algorithm 1}
\end{algorithm}

\textbf{Experiments for the shortest path and the maximum coverage problems.} %In the remainder of this section, we consider layered graphs with $h = 5$ intermediate layers and $r = 5$ nodes at each layer. 
As mentioned earlier, for both the SPP and the MCP we consider semi-bandit and bandit feedback scenarios described by the constraints (\ref{eq: data uncertainty spp}) and (\ref{eq: data uncertainty mcp}), respectively. As discussed in Section \ref{subsec: interval uncertainty}, the three-level problem~(\ref{three level formulation}) with the component-wise interval support constraints (\ref{eq: support constraints experiments}) and semi-bandit feedback can also be viewed as the SPP/MCP with interval uncertainty; 
in the case of bandit feedback, we need to solve the general MILP reformulation (\ref{mixed-integer programming reformulation}). 

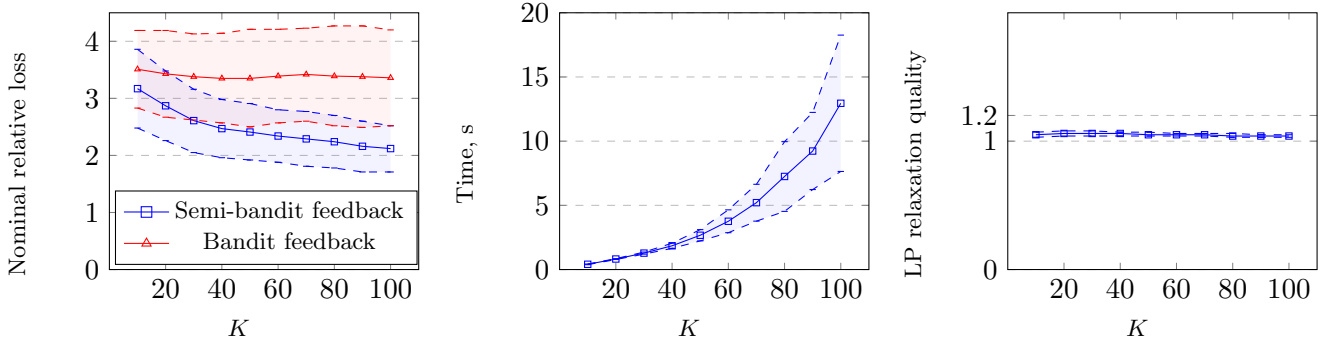
\begin{figure}
	\begin{center}
	\begin{subfigure}{0.3\textwidth}
	 \begin{tikzpicture}[scale = 0.6]
		\begin{axis}[
		xlabel=\footnotesize {$\qquad \quad K$},
		ylabel={ \footnotesize $\qquad$ Nominal relative loss},
		xmin=0, xmax=110,
		ymin=0, ymax=4.5,
		xtick={20,40,60,80,100},
		ytick={0,1,2,3,4},
		legend pos=south west,
		ymajorgrids=true,
		grid style=dashed,
		]
		\addplot[name path=f1,
		color=blue,
		mark=square,
		]
		coordinates {
			(10, 3.17)
			(20, 2.87)
			(30, 2.61)
			(40, 2.47)
			(50, 2.41)
			(60, 2.34)
			(70, 2.29)
			(80, 2.24)
			(90, 2.16)
			(100, 2.12)
			
		};
		
		\addplot[name path=f4,
		color=red,
		mark=triangle,
		]
		coordinates {
			(10, 3.51)
			(20, 3.43)
			(30, 3.38)
			(40, 3.35)
			(50, 3.35)
			(60, 3.39)
			(70, 3.42)
			(80, 3.39)
			(90, 3.38)
			(100, 3.36)
			
		};
		
		\addplot[name path=f2,
		color=blue,
		style=dashed,
		mark=-,
		]
		coordinates {
			(10, 2.48)
			(20, 2.26)
			(30, 2.05)
			(40, 1.96)
			(50, 1.92)
			(60, 1.88)
			(70, 1.81)
			(80, 1.78)
			(90, 1.71)
			(100, 1.71)	
		};
		\addplot[name path=f3,
		color=blue,
		style=dashed,
		mark=-,
		]
		coordinates {
		(10, 3.86)
		(20, 3.48)
		(30, 3.16)
		(40, 2.98)
		(50, 2.91)
		(60, 2.80)
	 (70, 2.77)
		(80, 2.70)
		(90, 2.60)
		(100, 2.52)
		};
		
		\addplot[name path=f5,
		style=dashed,
		color=red,
		mark=-,
		]
		coordinates {
			 (10, 2.83)
			 (20, 2.67)
			 (30, 2.62)
			 (40, 2.57)
			 (50, 2.50)
			 (60, 2.57)
			 (70, 2.60)
			 (80, 2.52)
			 (90, 2.49)
			 (100, 2.52)
			
		};
		
		\addplot[name path=f6,
		style=dashed,
		color=red,
		mark=-,
		]
		coordinates {
		 (10, 4.19)
		 (20, 4.19)
		 (30, 4.13)
		 (40, 4.14)
		 (50, 4.21)
		 (60, 4.21)
		 (70, 4.23)
		 (80, 4.27)
		 (90, 4.27)
		 (100, 4.20)
		};
		
		\addplot [
		thick,
		color=blue,
		fill=blue, 
		fill opacity=0.05
		]
		fill between[
		of=f3 and f2,
		soft clip={domain=5:100},
		];
		
		\addplot [
		thick,
		color=red,
		fill=red, 
		fill opacity=0.05
		]
		fill between[
		of=f6 and f5,
		soft clip={domain=5:100},
		];
		\legend{\footnotesize Semi-bandit feedback, \footnotesize Bandit feedback}
		\end{axis}
		\end{tikzpicture}
 \end{subfigure}	
	\hfill
		\begin{subfigure}{0.3\textwidth}
		\begin{tikzpicture}[scale = 0.6]
			\begin{axis}[
				xlabel=\footnotesize {$\qquad \quad K$},
				ylabel={\footnotesize $\qquad$ Time, s},
				xmin=0, xmax=110,
				ymin=0, ymax=20,
				xtick={20,40,60,80,100},
				ytick={0,5,10,15,20},
				legend pos=north west,
				ymajorgrids=true,
				grid style=dashed,
				]
				
				\addplot[name path=f1,
				color=blue,
				mark=square,
				]
				coordinates {
					(10, 0.41)
					(20, 0.82)
					(30, 1.28)
					(40, 1.85)
					(50, 2.67)
					(60, 3.76)
					(70, 5.21)
					(80, 7.25)
					(90, 9.23)
					(100, 12.95)
				};
				\addplot[name path=f2,
				color=blue,
				style=dashed,
				mark=-,
				]
				coordinates {
					(10, 0.38)
					(20, 0.78)
					(30, 1.20)
					(40, 1.65)
					(50, 2.23)
					(60, 2.87)
					(70, 3.77)
					(80, 4.53)
					(90, 6.23)
					(100, 7.64)
				};
				\addplot[name path=f3,
				color=blue,
				style=dashed,
				mark=-,
				]
				coordinates {
					(10, 0.43)
					(20, 0.85)
					(30, 1.36)
					(40, 2.05)
					(50, 3.11)
					(60, 4.65)
					(70, 6.64)
					(80, 9.96)
					(90, 12.24)
					(100, 18.26)
				};
				
				\addplot [
				thick,
				color=blue,
				fill=blue, 
				fill opacity=0.05
				]
				fill between[
				of=f3 and f2,
				soft clip={domain=10:100},
				];
			\end{axis}
		\end{tikzpicture}
	\end{subfigure}
	\hfill
	\begin{subfigure}{0.3\textwidth}
		\begin{tikzpicture}[scale = 0.6]
		\begin{axis}[
		xlabel=\footnotesize {$\qquad \quad K$},
		ylabel={\footnotesize $\qquad$ LP relaxation quality},
		xmin=0, xmax=110,
		ymin=0, ymax=2,
		xtick={20,40,60,80,100},
		ytick={0,1,1.2},
		legend pos=north west,
		ymajorgrids=true,
		grid style=dashed,
		]
		
		\addplot[name path=f1,
		color=blue,
		mark=square,
		]
		coordinates {
			(10, 1.05)
			(20, 1.06)
			(30, 1.06)
			(40, 1.06)
			(50, 1.05)
			(60, 1.05)
			(70, 1.05)
			(80, 1.04)
			(90, 1.04)
			(100, 1.04)
		};
		\addplot[name path=f2,
		color=blue,
		style=dashed,
		mark=-,
		]
		coordinates {
		(10, 1.03)
		(20, 1.04)
		(30, 1.04)
		(40, 1.04)
		(50, 1.04)
		(60, 1.04)
		(70, 1.04)
		(80, 1.03)
		(90, 1.03)
		(100, 1.03)
		};
		\addplot[name path=f3,
		color=blue,
		style=dashed,
		mark=-,
		]
		coordinates {
		(10, 1.07)
		(20, 1.08)
		(30, 1.08)
		(40, 1.07)
		(50, 1.07)
		(60, 1.06)
		(70, 1.06)
		(80, 1.06)
		(90, 1.05)
		(100, 1.05)
		};
		
		\addplot [
		thick,
		color=blue,
		fill=blue, 
		fill opacity=0.05
		]
		fill between[
		of=f3 and f2,
		soft clip={domain=10:100},
		];
		\end{axis}
		\end{tikzpicture}
	\end{subfigure}
	\caption{\footnotesize The SPP with $h = 11$, $r = 5$ and $K_{max} = 100$. For different types of feedback we report the average relative loss (\ref{eq: nominal relative loss}) with MADs as a function of the sample size, $K$, for $100$ random test instances. In the case of bandit feedback the average solution times and the average LP relaxation quality with MADs are also provided.}
	\label{fig: semi-bandit and bandit feedback SPP}
	\end{center}
\end{figure} 

In order to generate the set of paths $\mathcal{P}$ for the SPP, we employ a version of the upper confidence bound (UCB) algorithm for online COPs; see, e.g., \cite{Chen2013, Wang2017}. The pseudocode of the algorithm is given by Algorithm~\ref{alg: algorithm 1}. To put it briefly, in order to identify a path with the least expected cost we need (\textit{i}) to refine the expected cost of the best currently observed paths and (\textit{ii}) to examine the paths that are not sufficiently explored. In general, $\mathcal{P}$ can take on any form, but the aforementioned algorithm allows us to effectively extract new information from the observed~data. 

 It is worth noting that the decisions provided by (\ref{three level formulation}) do not account for an exploration stage and cannot be used directly in Algorithm~\ref{alg: algorithm 1}. Furthermore, the set of paths $\mathcal{P}$ in Algorithm \ref{alg: algorithm 1} is generated based on the assumption of semi-bandit feedback. We use \textit{the same} set of paths in the case of bandit feedback to provide a rather myopic comparison of these two settings. %More advanced online learning algorithms under the bandit feedback scenario are beyond the scope of the current study and can be found in \cite{Bubeck2012}. 
 
Finally, we use a similar algorithm in order to generate the set of decisions $\widetilde{\mathcal{P}}$ for the~MCP.

In the first experiment we explore the quality of distributionally robust decisions as a function of the sample size, $K$. In this regard, we sequentially increase the number of paths/decisions in $\mathcal{P}$ and $\widetilde{\mathcal{P}}$, respectively, and report the average relative loss (\ref{eq: nominal relative loss}) with MADs for 100 random test instances. In the case of semi-bandit feedback, we use Theorem \ref{theorem 2}; in the case of bandit feedback, we solve the MILP reformulation (\ref{mixed-integer programming reformulation}) providing the average solution time and the average LP relaxation quality (i.e., the ratio of optimal objective function values of (\ref{mixed-integer programming reformulation}) and its linear programming relaxation).

We set $r = 5$, $h = 11$, $K_{max} = 100$ and $\varepsilon = \sqrt{\frac{K_{max}}{K}}$ for the SPP. Also, we set $n_1 = n_2 = 50$, $\tilde{h} = 5$, $K_{max} = 50$ and $\varepsilon = \sqrt{\frac{K_{max}}{K}}$ for the MCP. Subsets $\mathcal{A}_i \subseteq \mathcal{A}$, $i \in \{1, \ldots, n_2\}$, for the latter problem are assumed to have a fixed cardinality, $5$, and are generated uniformly at random. We also note that since the MCP is a maximization problem, we set $\ell(\mathbf{x}, \mathbf{c}) = -\mathbf{c}^\top \mathbf{x}$ in Theorem \ref{theorem 1} and reformulate Theorem \ref{theorem 2} symmetrically with respect to the lower bounds. Consequently, both the nominal relative loss (\ref{eq: nominal relative loss}) and the LP relaxation quality for the MCP do not exceed~$1$. The results for the SPP and the MCP are reported in Figures \ref{fig: semi-bandit and bandit feedback SPP} and \ref{fig: semi-bandit and bandit feedback MCP}, respectively. 

 The observations can be summarized as follows:
\begin{itemize}
	\item For both problems and both types of information feedback the nominal relative loss (\ref{eq: nominal relative loss}) approaches to $1$ with the increase of the number of samples, $K$. However, the convergence rate is less significant in the case of bandit feedback. This observation can be attributed to the fact that as $K$ increases, the decision-maker gains more distributional information. However, the quality of this information is lower in the case of bandit feedback. 
	\item From Figure \ref{fig: semi-bandit and bandit feedback SPP} we observe that the average solution times for the SPP increase almost linearly for reasonably small $K$, and the average LP relaxation quality is close to $1$. In this regard, we note the number of non-zero components for any $\mathbf{x} \in X^{(SPP)}$ is sufficiently small compared to the number of arcs in the graph and, hence, it can often be the case that the observed paths in $\mathcal{P}$ do not overlap with each other. This observation along with Theorem \ref{theorem 3} provide some intuition behind the fact that the SPP with bandit feedback can be solved rather effectively for sufficiently small values of $K$. 
	\item On the other hand, as shown in Figure \ref{fig: semi-bandit and bandit feedback MCP}, the average solution times for the MCP are inversely proportional to the LP relaxation quality; this trend is more pronounced than the expected increase of solution times in $K$. In this regard, we note that the number of non-zero components for any $\mathbf{x} \in X^{MCP}$ is relatively large in terms of~$n_1$. Hence, the conditions of Theorem \ref{theorem 3} are not satisfied and the solution times increase substantially even for small values of $K$. 
	
\begin{figure}
	\begin{center}
			\begin{subfigure}{0.3\textwidth}
			\begin{tikzpicture}[scale = 0.6]
				\begin{axis}[
					xlabel=\footnotesize {$\qquad \quad K$},
					ylabel={ \footnotesize $\qquad$ Nominal relative loss},
					xmin=0, xmax=60,
					ymin=0, ymax=1.5,
					xtick={10,20,30,40,50},
					ytick={0,0.5,0.75,1},
					legend pos=south west,
					ymajorgrids=true,
					grid style=dashed,
					]
					\addplot[name path=f1,
					color=blue,
					mark=square,
					]
					coordinates {
						(5, 0.92)
						(10, 0.96)
						(15, 0.97)
						(20, 0.98)
						(25, 0.98)
						(30, 0.98)
						(35, 0.98)
						(40, 0.98)
						(45, 0.98)
						(50, 0.98)
					};
					
					\addplot[name path=f4,
					color=red,
					mark=triangle,
					]
					coordinates {
						(5, 0.80)
						(10, 0.86)
						(15, 0.91)
						(20, 0.93)
						(25, 0.95)
						(30, 0.96)
						(35, 0.96)
						(40, 0.96)
						(45, 0.97)
						(50, 0.97)
					};
					
					\addplot[name path=f2,
					color=blue,
					style=dashed,
					mark=-,
					]
					coordinates {
						(5, 0.89)
						(10, 0.93)
						(15, 0.95)
						(20, 0.96)
						(25, 0.96)
						(30, 0.96)
						(35, 0.97)
						(40, 0.97)
						(45, 0.97)
						(50, 0.97)
					};
					\addplot[name path=f3,
					color=blue,
					style=dashed,
					mark=-,
					]
					coordinates {
						(5, 0.95)
						(10, 0.99)
						(15, 0.99)
						(20, 1.00)
						(25, 1.00)
						(30, 1.00)
						(35, 1.00)
						(40, 1.00)
						(45, 1.00)
						(50, 1.00)
					};
					
					\addplot[name path=f5,
					style=dashed,
					color=red,
					mark=-,
					]
					coordinates {
						(5, 0.74)
						(10, 0.81)
						(15, 0.87)
						(20, 0.90)
						(25, 0.92)
						(30, 0.93)
						(35, 0.94)
						(40, 0.94)
						(45, 0.95)
						(50, 0.95)
					};
					
					\addplot[name path=f6,
					style=dashed,
					color=red,
					mark=-,
					]
					coordinates {
						(5, 0.86)
						(10, 0.91)
						(15, 0.95)
						(20, 0.97)
						(25, 0.98)
						(30, 0.98)
						(35, 0.98)
						(40, 0.99)
						(45, 0.99)
						(50, 0.99)
					};
					
					\addplot [
					thick,
					color=blue,
					fill=blue, 
					fill opacity=0.05
					]
					fill between[
					of=f3 and f2,
					soft clip={domain=5:100},
					];
					
					\addplot [
					thick,
					color=red,
					fill=red, 
					fill opacity=0.05
					]
					fill between[
					of=f6 and f5,
					soft clip={domain=5:100},
					];
					\legend{\footnotesize Semi-bandit feedback, \footnotesize Bandit feedback}
				\end{axis}
			\end{tikzpicture}
		\end{subfigure}
		\hfill
			\begin{subfigure}{0.3\textwidth}
			\begin{tikzpicture}[scale = 0.6]
				\begin{axis}[
					xlabel=\footnotesize {$\qquad \quad K$},
					ylabel={\footnotesize $\qquad$ Time, s},
					xmin=0, xmax=60,
					ymin=0, ymax=85,
					xtick={10,20,30,40,50},
					ytick={0,10,20,30,40,80},
					legend pos=north west,
					ymajorgrids=true,
					grid style=dashed,
					]
					
					\addplot[name path=f1,
					color=blue,
					mark=square,
					]
					coordinates {
						(5, 17.12)
						(10, 38.20)
						(15, 20.04)
						(20, 11.05)
						(25, 10.27)
						(30, 6.99)
						(35, 6.21)
						(40, 6.28)
						(45, 5.94)
						(50, 5.34)
					};
					\addplot[name path=f2,
					color=blue,
					style=dashed,
					mark=-,
					]
					coordinates {
						(5, -5.93)
						(10, -10.24)
						(15, -2.45)
						(20, 0.22)
						(25, -0.54)
						(30, 0.22)
						(35, 0.17)
						(40, 0.12)
						(45, 0.16)
						(50, 0.80)
					};
					\addplot[name path=f3,
					color=blue,
					style=dashed,
					mark=-,
					]
					coordinates {
						(5, 40.18)
						(10, 86.64)
						(15, 42.53)
						(20, 21.87)
						(25, 21.08)
						(30, 13.76)
						(35, 12.24)
						(40, 12.44)
						(45, 11.73)
						(50, 9.88)
					};
					
					\addplot [
					thick,
					color=blue,
					fill=blue, 
					fill opacity=0.05
					]
					fill between[
					of=f3 and f2,
					soft clip={domain=5:50},
					];
				\end{axis}
			\end{tikzpicture}
		\end{subfigure}
		\hfill
		\begin{subfigure}{0.3\textwidth}
			\begin{tikzpicture}[scale = 0.6]
				\begin{axis}[
					xlabel=\footnotesize {$\qquad \quad K$},
					ylabel={\footnotesize $\qquad$ LP relaxation quality},
					xmin=0, xmax=60,
					ymin=0, ymax=2,
					xtick={10,20,30,40,50},
					ytick={0, 0.4, 0.6, 0.8, 1},
					legend pos=north west,
					ymajorgrids=true,
					grid style=dashed,
					]
					
					\addplot[name path=f1,
					color=blue,
					mark=square,
					]
					coordinates {
						(5, 0.41)
						(10, 0.36)
						(15, 0.43)
						(20, 0.51)
						(25, 0.58)
						(30, 0.63)
						(35, 0.67)
						(40, 0.70)
						(45, 0.72)
						(50, 0.74)
					};
					\addplot[name path=f2,
					color=blue,
					style=dashed,
					mark=-,
					]
					coordinates {
						(5, 0.21)
						(10, 0.22)
						(15, 0.32)
						(20, 0.41)
						(25, 0.48)
						(30, 0.53)
						(35, 0.58)
						(40, 0.60)
						(45, 0.63)
						(50, 0.65)
					};
					\addplot[name path=f3,
					color=blue,
					style=dashed,
					mark=-,
					]
					coordinates {
						(5, 0.62)
						(10, 0.51)
						(15, 0.55)
						(20, 0.62)
						(25, 0.68)
						(30, 0.73)
						(35, 0.77)
						(40, 0.79)
						(45, 0.81)
						(50, 0.82)
					};
					
					\addplot [
					thick,
					color=blue,
					fill=blue, 
					fill opacity=0.05
					]
					fill between[
					of=f3 and f2,
					soft clip={domain=10:100},
					];
				\end{axis}
			\end{tikzpicture}
		\end{subfigure}
		\caption{\footnotesize The MCP with $n_1 = n_2 = 50$, $\tilde{h} = 5$ and $K_{max} = 50$. For different types of feedback we report the average relative loss (\ref{eq: nominal relative loss}) with MADs as a function of the sample size, $K$, for $100$ random test instances. In the case of bandit feedback the average solution times and the average LP relaxation quality with MADs are also provided. }
		\label{fig: semi-bandit and bandit feedback MCP}
		\end{center}
	\end{figure}
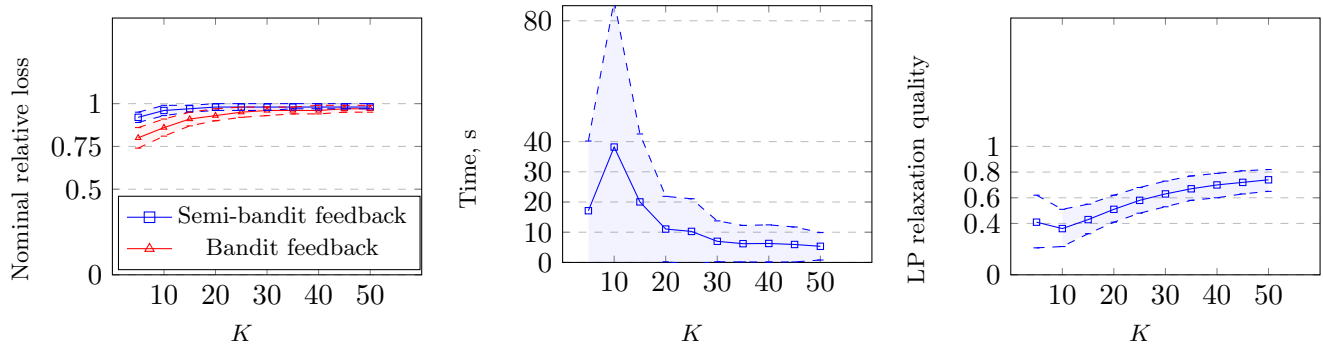
	
\item The last two observations also imply that the rate of convergence for the nominal relative loss~(\ref{eq: nominal relative loss}) is faster for the MCP than for the SPP. In fact, the more non-zero components are contained in $\mathbf{x} \in \{0, 1\}^n$, the more qualitative information is typically observed from the costs associated with this decision. \end{itemize}

\begin{figure}
	\begin{center}
		\begin{subfigure}{0.3\textwidth}
			\begin{tikzpicture}[scale = 0.6]
			\begin{axis}[
			xlabel=\footnotesize {$\qquad \quad h$},
			ylabel={ \footnotesize $\qquad$ Nominal relative loss},
			xmin=0, xmax=18,
			ymin=0, ymax=5,
			xtick={1,5,9,13,17},
			ytick={0,1,2,3,4},
			legend pos=south west,
			ymajorgrids=true,
			grid style=dashed,
			]
			\addplot[name path=f1,
			color=blue,
			mark=square,
			]
			coordinates {
				(5, 2.05)
				(7, 2.61)
				(9, 2.88)
				(11, 3.17)
				(13, 3.42)
				(15, 3.52)
				(17, 3.68)	
			};
			\addplot[name path=f2,
			color=blue,
			style=dashed,
			mark=-,
			]
			coordinates {
				(5, 1.46)
				(7, 1.87)
				(9, 2.12)
				(11, 2.50)
				(13, 2.69)
				(15, 2.89)
				(17, 3.04)	
			};
			\addplot[name path=f3,
			color=blue,
			style=dashed,
			mark=-,
			]
			coordinates {
				(5, 2.64)
				(7, 3.34)
				(9, 3.64)
				(11, 3.85)
				(13, 4.14)
				(15, 4.16)
				(17, 4.32)
			};
		
			\addplot [
			thick,
			color=blue,
			fill=blue, 
			fill opacity=0.05
			]
			fill between[
			of=f3 and f2,
			soft clip={domain=5:17},
			];
			\end{axis}
			\end{tikzpicture}
		\end{subfigure}	
		\hfill
		\begin{subfigure}{0.3\textwidth}
			\begin{tikzpicture}[scale = 0.6]
			\begin{axis}[
			xlabel=\footnotesize {$\qquad \quad h$},
			ylabel={\footnotesize $\qquad$ Time, s},
			xmin=0, xmax=18,
			ymin=0, ymax=210,
		 xtick={1,5,9,13,17},
			ytick={0,50,100,150,200},
			legend pos=north west,
			ymajorgrids=true,
			grid style=dashed,
			]
			
			\addplot[name path=f1,
			color=blue,
			mark=square,
			]
			coordinates {
				(5, 0.32)
				(7, 0.75)
				(9, 1.42)
				(11, 2.71)
				(13, 7.58)
				(15, 21.55)
				(17, 111.37)
			}; 
			\addplot[name path=f2,
			color=blue,
			style=dashed,
			mark=-,
			]
			coordinates {
				(5,0.28)
				(7,0.68)
				(9,1.30)
				(11,2.20)
				(13,5.33)
				(15,8.07)
				(17,29.58)
			};
			\addplot[name path=f3,
			color=blue,
			style=dashed,
			mark=-,
			]
			coordinates {
				(5,0.36)
				(7,0.83)
				(9,1.54)
				(11,3.22)
				(13,9.82)
				(15,35.04)
				(17,192.17)
			};
			
			\addplot [
			thick,
			color=blue,
			fill=blue, 
			fill opacity=0.05
			]
			fill between[
			of=f3 and f2,
			soft clip={domain=5:17},
			];
			\end{axis}
			\end{tikzpicture}
		\end{subfigure}
		\hfill
		\begin{subfigure}{0.3\textwidth}
			\begin{tikzpicture}[scale = 0.6]
			\begin{axis}[
			xlabel=\footnotesize {$\qquad \quad h$},
			ylabel={\footnotesize $\qquad$ LP relaxation quality},
			xmin=0, xmax=18,
			ymin=0, ymax=2,
	 	xtick={1,5,9,13,17},
			ytick={0,1,1.2},
			legend pos=north west,
			ymajorgrids=true,
			grid style=dashed,
			]
			
			\addplot[name path=f1,
			color=blue,
			mark=square,
			]
			coordinates {
				(5,1.04)
				(7,1.05)
				(9,1.05)
				(11,1.05)
				(13,1.05)
				(15,1.05)
				(17,1.05)
			};
			\addplot[name path=f2,
			color=blue,
			style=dashed,
			mark=-,
			]
			coordinates {
				(5,1.01)
				(7,1.02)
				(9,1.04)
				(11,1.04)
				(13,1.04)
				(15,1.04)
				(17,1.04)
			};
			\addplot[name path=f3,
			color=blue,
			style=dashed,
			mark=-,
			]
			coordinates {
				(5,1.07)
				(7,1.07)
				(9,1.07)
				(11,1.06)
				(13,1.06)
				(15,1.06)
				(17,1.05)
			};
			
			\addplot [
			thick,
			color=blue,
			fill=blue, 
			fill opacity=0.05
			]
			fill between[
			of=f3 and f2,
			soft clip={domain=5:17},
			];
			\end{axis}
			\end{tikzpicture}
		\end{subfigure}
		\caption{\footnotesize The SPP with $r = 5$, $K = 50$ and $\varepsilon_K = \frac{\sqrt{2}}{11}h$. We report the average relative loss (\ref{eq: nominal relative loss}), the average solution times and the average LP relaxation quality with MADs as a function of $h$ for $100$ random test instances.}
		\label{fig: dependence on the size SPP}
	\end{center}
\end{figure}
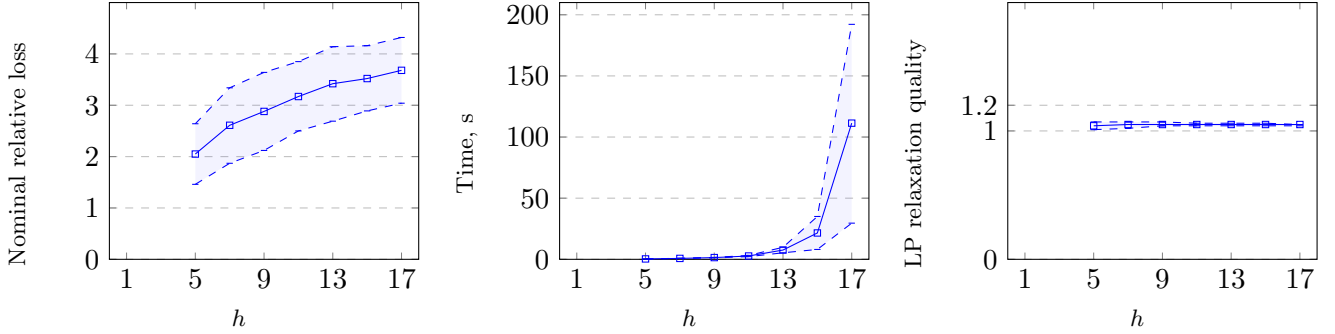

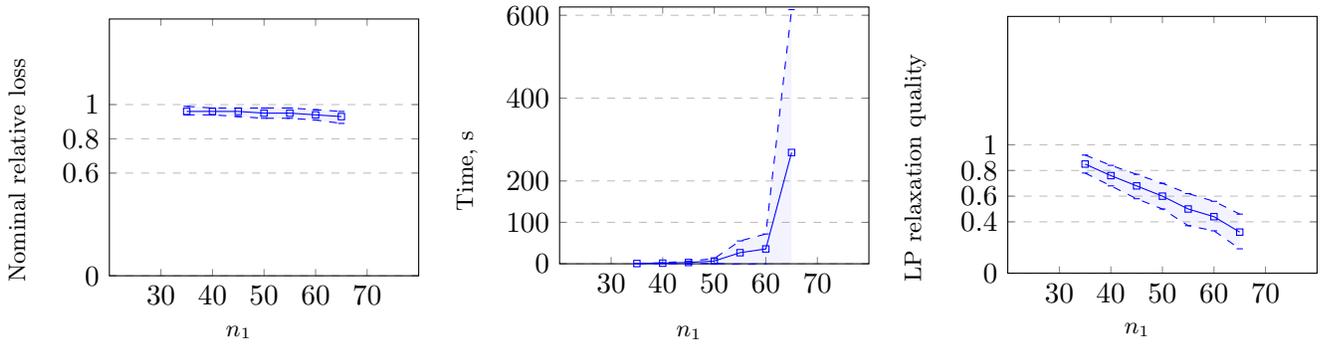
\begin{figure}
	\begin{center}
		\begin{subfigure}{0.3\textwidth}
			\begin{tikzpicture}[scale = 0.6]
				\begin{axis}[
					xlabel=\footnotesize {$\qquad \quad n_1$},
					ylabel={ \footnotesize $\qquad$ Nominal relative loss},
					xmin=20, xmax=80,
					ymin=0, ymax=1.5,
					xtick={30,40,50,60,70},
					ytick={0,0.6,0.8,1},
					legend pos=south west,
					ymajorgrids=true,
					grid style=dashed,
					]
					\addplot[name path=f1,
					color=blue,
					mark=square,
					]
					coordinates {
						(35, 0.96)
						(40, 0.96)
						(45, 0.96)
						(50, 0.95)
						(55, 0.95)
						(60, 0.94)
						(65, 0.93)	
					};
					\addplot[name path=f2,
					color=blue,
					style=dashed,
					mark=-,
					]
					coordinates {
						(35, 0.94)
						(40, 0.94)
						(45, 0.93)
						(50, 0.92)
						(55, 0.92)
						(60, 0.91)
						(65, 0.89)	
					};
					\addplot[name path=f3,
					color=blue,
					style=dashed,
					mark=-,
					]
					coordinates {
						(35, 0.99)
						(40, 0.98)
						(45, 0.98)
						(50, 0.98)
						(55, 0.98)
						(60, 0.97)
						(65, 0.96)
					};
					
					\addplot [
					thick,
					color=blue,
					fill=blue, 
					fill opacity=0.05
					]
					fill between[
					of=f3 and f2,
					soft clip={domain=35:65},
					];
				\end{axis}
			\end{tikzpicture}
		\end{subfigure}	
		\hfill
		\begin{subfigure}{0.3\textwidth}
			\begin{tikzpicture}[scale = 0.6]
				\begin{axis}[
					xlabel=\footnotesize {$\qquad \quad n_1$},
					ylabel={\footnotesize $\qquad$ Time, s},
					xmin=20, xmax=80,
					ymin=0, ymax=620,
					xtick={30,40,50,60,70},
					ytick={0,100,200,400,600},
					legend pos=north west,
					ymajorgrids=true,
					grid style=dashed,
					]
					
					\addplot[name path=f1,
					color=blue,
					mark=square,
					]
					coordinates {
						(35, 0.58)
						(40, 1.43)
						(45, 3.15)
						(50, 6.61)
						(55, 26.71)
						(60, 35.79)
						(65, 268.68)
					};
					\addplot[name path=f2,
					color=blue,
					style=dashed,
					mark=-,
					]
					coordinates {
						(35, 0.19)
						(40, 0.26)
						(45, 0.16)
						(50, 0.59)
						(55, -1.75)
						(60, -0.23)
						(65, -76.26)
					};
					\addplot[name path=f3,
					color=blue,
					style=dashed,
					mark=-,
					]
					coordinates {
						(35, 0.96)
						(40, 2.61)
						(45, 6.14)
						(50, 12.63)
						(55, 55.17)
						(60, 71.82)
						(65, 613.62)
					};
					
					\addplot [
					thick,
					color=blue,
					fill=blue, 
					fill opacity=0.05
					]
					fill between[
					of=f3 and f2,
					soft clip={domain=35:65},
					];
				\end{axis}
			\end{tikzpicture}
		\end{subfigure}
		\hfill
		\begin{subfigure}{0.3\textwidth}
			\begin{tikzpicture}[scale = 0.6]
				\begin{axis}[
					xlabel=\footnotesize {$\qquad \quad n_1$},
					ylabel={\footnotesize $\qquad$ LP relaxation quality},
					xmin=20, xmax=80,
					ymin=0, ymax=2,
					xtick={30,40,50,60,70},
					ytick={0,0.4,0.6,0.8,1},
					legend pos=north west,
					ymajorgrids=true,
					grid style=dashed,
					]
					
					\addplot[name path=f1,
					color=blue,
					mark=square,
					]
					coordinates {
						(35, 0.85)
						(40, 0.76)
						(45, 0.68)
						(50, 0.60)
						(55, 0.50)
						(60, 0.44)
						(65, 0.32)
					};
					\addplot[name path=f2,
					color=blue,
					style=dashed,
					mark=-,
					]
					coordinates {
						(35, 0.78)
						(40, 0.68)
						(45, 0.58)
						(50, 0.50)
						(55, 0.37)
						(60, 0.33)
						(65, 0.19)
					};
					\addplot[name path=f3,
					color=blue,
					style=dashed,
					mark=-,
					]
					coordinates {
						(35, 0.92)
						(40, 0.84)
						(45, 0.77)
						(50, 0.70)
						(55, 0.62)
						(60, 0.56)
						(65, 0.46)
					};
					
					\addplot [
					thick,
					color=blue,
					fill=blue, 
					fill opacity=0.05
					]
					fill between[
					of=f3 and f2,
					soft clip={domain=35:65},
					];
				\end{axis}
			\end{tikzpicture}
		\end{subfigure}
		\caption{\footnotesize The MCP with $n_2 = 50$, $\tilde{h} = 5$, $K = 25$ and $\varepsilon_K = \frac{\sqrt{2}}{50}n_1$. We report the average relative loss (\ref{eq: nominal relative loss}), the average solution times and the average LP relaxation quality with MADs as a function of $n_1$ for $100$ random test instances. }
		\label{fig: dependence on the size MCP}
	\end{center}
\end{figure}

In our last experiment, we assume the bandit feedback scenario and consider the nominal relative loss (\ref{eq: nominal relative loss}) as a function of the size of each problem. Specifically, we set $r = 5$ and $K = 50$ for the~SPP and explore how the quality of distributionally robust decisions and solution times for the MILP reformulation \ref{mixed-integer programming reformulation}) scale in the length of each path, $h$. Analogously, we set $n_2 = 50$, $\tilde{h} = 5$, $K = 25$ for the MCP and explore the dependence on the number of items, $n_1$ (the subsets $\mathcal{A}_i$, $i \in \{1, \ldots, n_2\}$, are selected precisely as in the previous experiment). Finally, according to equation (\ref{eq: Hoeffding}), we need to select $\varepsilon_K$ proportional to $\max_{\mathbf{x} \in X} \Vert \mathbf{x} \Vert_1$ and, hence, we set $\varepsilon_K = \gamma' h$ for the SPP and $\varepsilon_K = \gamma'' n_1$ for the MCP, where $\gamma' \in \mathbb{R}_{++}$ and $\gamma'' \in \mathbb{R}_{++}$ are selected consistently with the previous experiment. The results are reported in Figures \ref{fig: dependence on the size SPP} and \ref{fig: dependence on the size MCP}. 

We make the following observations:
\begin{itemize}
	\item The degree of misspecification of the optimal solution of (\ref{stochastic programming problem}) increases with the increase of problem's size for both the SPP and the MCP. This observation is rather intuitive as the sample size, $K$, is fixed. Hence, the smaller the problem's size is, the more effectively the observed costs can be used in reconstruction of the optimal solution. 
	\item Similar to the previous experiment, from Figure \ref{fig: dependence on the size SPP} we observe that the LP relaxation quality for the SPP is close to $1$ and the problem can be solved rather effectively even for graphs of a medium size. 
	
	\item On the other hand, as depicted in Figure \ref{fig: dependence on the size MCP}, solution times for the MCP increase faster than for the SPP; we recall that, in contrast to the SPP, the nominal version of the MCP is known to be $NP$-hard. 
\end{itemize}
Finally, from Figures \ref{fig: semi-bandit and bandit feedback MCP} and \ref{fig: dependence on the size MCP} we observe that the LP relaxation quality for the MCP non-trivially depends on the number of samples, $K$, and tends to decrease in the number of items, $n_1$. We suppose that these dependencies are stipulated by a particular structure of the MCP and, hence, no valuable conclusions can be made for the general case. 

%As expected, both the nominal relative loss (\ref{eq: nominal relative loss}) and the solution times increase in the size of the graph. However, it is worth noting that the LP relaxation quality remains sufficiently small, irrespective of the size of the graph. Hence, in most cases the three-level shortest path problem (\ref{three level formulation}) with bandit feedback can be solved reasonably fast even for graphs of a medium size (however, there are a few test instances, for which solution times are sufficiently large). 

\section{Conclusions} \label{sec: conclusions}
In this paper we examine a class of linear mixed-integer programming (MILP) problems, where the distribution of the cost vector is subject to uncertainty. In contrast to most of the related studies, this distribution can only be observed through a finite data set, which itself is subject to uncertainty. The overall problem is formulated as a three-level distributionally robust optimization (DRO) problem, where the decision-maker aims to minimize its expected loss under the worst-case realization of both the data and the data-generating~distribution.

Our approach to modeling data uncertainty exploits sample-wise linear constraints that are tailored to the optimization problem's structure and the data collection process. From a practical perspective, it allows us to effectively capture various particular forms of data uncertainty, such as noise, misspecification and incomplete information feedback. In particular, the latter form of data uncertainty is motivated in the context of related online combinatorial optimization problem settings. %Additionally, to describe the uncertainty in the distribution for a fixed realization of data, we use a type-1 Wasserstein ball centered at the empirical distribution of the data. 
%The key contribution of this study is that we propose a novel approach to modeling data uncertainty in the context of linear mixed-integer DRO problems. In this regard, whenever the cost vector is subject to linear support constraints, we describe each random sample by another set of \textit{linear constraints} tailored to the optimization problem's structure and/or data collection process. From a practical perspective, it is shown that various particular forms of data uncertainty, such as interval uncertainty, misspecification and incomplete information feedback, can be described in a rather straightforward way by using linear constraints with respect to the data. 

\looseness-1 We demonstrate that the proposed three-level DRO problem with an $l_1$-norm Wasserstein ambiguity set and a biaffine loss function can be reformulated as a single-level MILP problem. Moreover, for two particular forms of data uncertainty, interval uncertainty and ``bandit'' feedback, two specific reformulations are provided. These reformulations are limited to combinatorial optimization problems of a predefined structure but can be solved rather effectively using the underlying deterministic~problem. 

Finally, the obtained theoretical results are used in our computational study, where the three-level optimization model is applied to several classes of stochastic combinatorial optimization problems. First, we observe that well-organized incomplete/partially observable data allows to improve the model's out-of-sample performance. Secondly, the MILP reformulation with ``bandit'' feedback can be solved reasonably fast, if the sample size is relatively small and the observed decisions are sparse in terms of the problem's dimension. 

With respect to future research directions, it would be interesting to explore another forms of data-driven ambiguity sets in the context of DRO problems with incomplete data. Additionally, although our approach is focused on a single-level MILP reformulation of the three-level problem, one may suggest using existing algorithms to min-max-min robust MIP problems; see, e.g., \cite{Subramanyam2020}. It seems that these algorithms may allow to consider more general loss functions, e.g., piecewise linear concave functions in the uncertain problem parameters, providing however that the three-level problem can only be solved approximately. 

\textbf{Acknowledgments.} The article was prepared within the framework of the Basic Research Program at the National Research University Higher School of Economics (Sections 1-2). The part of the work carried out at HSE university was funded by RSF grant №22-11-00073~(Sections 3-5). The research for this paper is not based on or part of a cooperation agreement between the mentioned institution and the University of Zurich. 

\textbf{Conflict of interest:} None.

\textbf{Data availability statement:} The author confirms that all data generated or analysed during this study are included in this published article.

\bibliographystyle{apa}
\bibliography{bibliography}

\end{document}